\documentclass[12pt]{article}
\usepackage{amsmath, amssymb, latexsym, amscd, amsthm,amsfonts,amstext}
\usepackage[mathscr]{eucal}
\usepackage{graphicx}
\usepackage{subfig}

\newtheorem{theorem}{Theorem}[section]
\newtheorem{lemma}[theorem]{Lemma}

\newtheorem{remark}{Remark}[section]
\numberwithin{equation}{section}
\makeatletter
\renewcommand\tableofcontents{    \@starttoc{toc}}
\makeatother
\DeclareMathOperator{\Imag}{Im}
\DeclareMathOperator{\Real}{Re}

\def\bbH{\mathbb{H}}

\begin{document}

\title{Solving a 1-D inverse medium scattering problem using a new
multi-frequency globally strictly convex objective functional}
\author{Nguyen T. Th\`anh$^{1*}$ and Michael V. Klibanov$^2$ \\
\\
$^1$Department of Mathematics, Rowan University, \\
201 Mullica Hill Rd, Glassboro, NJ 08028, USA. \\
Email: \texttt{nguyent@rowan.edu} \\
$^2$Department of Mathematics and Statistics, \\
University of North Carolina at Charlotte\\
9201 University City Blvd, Charlotte, NC, 28223, USA. \\
Email: \texttt{mklibanv@uncc.edu}\\
$^*$Corresponding author. 
}

\date{}
\maketitle

\begin{abstract}
We propose a new approach to constructing globally strictly convex objective
functional in a 1-D inverse medium scattering problem using multi-frequency
backscattering data. The global convexity of the proposed objective functional is proved
using a Carleman estimate. Due to its convexity, no good first guess is
required in minimizing this objective functional. We also prove the global
convergence of the gradient projection algorithm and derive an error
estimate for the reconstructed coefficient. Numerical results show
reasonable reconstruction accuracy for simulated data.
\end{abstract}


%
%
%
%
%
%
%
%

\section{Introduction}

\label{sec:int}

One of the most popular techniques used for the purpose of detection of
buried objects is the Ground Penetrating Radar (GPR). Exploiting the energy
of backscattering electromagnetic pulses measured on the ground, the GPR
allows for mapping underground structures. The radar community mainly uses
migration-type imaging methods to obtain geometrical information such as the
shapes, the sizes, and the locations of the targets, see, e.g., \cite%
{C-C:2006, Daniels:2004,Ito:IP12,Ito:IP13, Li:IP15, Li:JCP14,Li:IP07,  Soumekh:1999,Yilmaz:1987}. However, these methods cannot
determine physical characteristics of buried objects. Therefore, additional
information about the objects' physical properties, such as the dielectric
permittivity and the magnetic permeability, may be helpful for their
identification.

The problem of determining these parameters can be formulated as a
coefficient identification problem for the wave equation. In the scattering
theory, this is also called an inverse medium scattering problem. This problem
has been extensively investigated, see e.g. \cite{C-K:2013} and the
references therein. Several methods have been proposed for solving it. One
of the earliest approach is the Born approximation which is effective at low
frequencies, see e.g., \cite{Bleistein:1984}. For gradient-based and
Newton-type methods, we refer the reader to \cite{Chavent:2009, Gon:IP17, Hohage:IP2001, Lakhal:IP10, NW:IP1995,Riz:IP17}
and the references therein. For decomposition methods, see e.g., \cite%
{C-K:2013, CM:QJMAM1988}.

%
%
%

In this work, we consider an inverse medium scattering problem in one
dimension using backscattering data generated by a single source position at multiple
frequencies. The model is described by the following equation: 
\begin{equation}
u''\left( x\right) +k^{2}c(x)u(x)=-\delta (x-x^{0}),\
x\in \mathbb{R},  \label{1.0}
\end{equation}%
where $k$ is the wavenumber, $c(x)$ represents the dielectric constant of
the medium in which the wave, originated by the point source at $x^{0},$
travels. The purpose of the coefficient inverse problem (CIP) under
consideration is to determine the coefficient $c(x)$ from measurements of $%
u(x,k)$ at a single location associated with multiple frequencies. One of
more specific applications is in the identification of mine-like targets.\
In this instance we refer to works with experimental data measured in the
field by a forward looking radar of US Army Research Laboratory \cite{KK:CMA2018,KKNS:SIAP2017,  KKNS:IP18}.

Using the multi-frequency data in inverse scattering problems has been
reported to be efficient. There are
different ways of using multi-frequency data. One approach, known as
frequency-hopping algorithms, uses the reconstruction at a lower frequency
as an initial guess for the reconstruction at a higher frequency. Several results have been
reported, see e.g., \cite{BLT:IP2015,B-B-P-S-S:JEWA2000, Chen:IP1997,
C-L:IEEE1995, S-T:IPI2012, ST:ESAIMM2NA2014, T-B-L-H:IP2001,
T-B-L-H:IEEETGRS2001}. Another approach is to use non-iterative
sampling-type methods, see, e.g.~\cite{Griesmaier:IP2011, GS:SIIM2017,
G-C-B:IP2010, Potthast:IP2010}. The third type of methods, based on the
construction of globally strictly convex objective functionals or a globally
convergent iterative process, has been reported recently, see, e.g.~\cite%
{KK:CMA2018,KKNS:SIAP2017, KKNS:IP18, KNNL:IPI2018,
KKNNN:ANM2017,NKNF:JIIP2018}.

In this paper, we continue our research on the third approach. The key step
of this method is the construction of an objective functional which contains
a Carleman Weight Function (CWF). The key property of this functional is
that it is strictly convex on any given set in an appropriate function space
if the parameter of the CWF is chosen large enough. This makes the method
converge globally, which is unlike conventional optimization-based
approaches which are usually locally convergent.

The idea of this type of methods was investigated earlier in \cite%
{KT:SIAP2015} for a similar problem in time domain. Then it was developed for multi-frequency measurements
in \cite{KK:CMA2018,KKNS:SIAP2017, KKNS:IP18}. In the time-domain
problem, the forward scattering model is described by the following Cauchy
problem: 
\begin{eqnarray}
&&c\left( x\right) u_{tt}=u_{xx}+\delta (x-x^{0},t),\ \left( x,t\right) \in 
\mathbb{R}\times \left( 0,\infty \right) ,  \label{6.1} \\
&&u\left( x,0\right) =0,u_{t}\left( x,0\right) =0.  \label{6.2}
\end{eqnarray}%
To reconstruct the coefficient $c(x)$,  in \cite{KT:SIAP2015} we established a globally strictly convex objective functional in the Laplace transform
domain. More precisely, let $s\geq \underline{s}>0$ be the Laplace transform
parameter, which is usually referred to as the \textit{pseudo frequency}.
The Laplace transform $\tilde{u}(x,s)$ of $u(x,t)$ satisfies the following
equation: 
\begin{equation}
\tilde{u}_{xx}(x,s)-s^{2}c(x)\tilde{u}(x,s)=-\delta (x-x^{0}),\ x\in \mathbb{%
R},  \label{1.01}
\end{equation}%
It can be proved that $\tilde{u}>0$ for $s$ large enough. By defining new
functions 
\begin{equation}
v(x,s):=\ln (\tilde{u})/s^{2};\quad q(x,s):=\dfrac{\partial v(x,s)}{\partial
s},  \label{eq1.5}
\end{equation}%
we obtain a nonlinear integro-differential equation for $q$. This equation
does not contain the unknown coefficient $c(x)$. However, $c(x)$ can be
calculated if $q(x,s)$ is known. The problem of finding $q(x,s)$ is then
converted to a minimization problem in which the objective functional is
globally strictly convex. The methods in \cite{KK:CMA2018,KKNS:SIAP2017,
KKNS:IP18} are similar, except that the Laplace transform is
not needed since the frequency domain problem can be treated as obtained from the time domain problem by the Fourier transform. There is an advantage
of the frequency-domain approach compared to the time-domain one is that in
the time-domain model, only signals which arrive at the receiver
early are usable in the inverse problems. This is because the kernel of the
Laplace transform decays exponentially in time. As a result, information
contained in later signals is diminished after the Laplace transform.
Consequently, the reconstruction accuracy is good only near the location of
measurement. This is not the case for the frequency-domain data. 

However, in the frequency-domain approach of \cite%
{KK:CMA2018,KKNS:SIAP2017,KKNS:IP18}, the solution of \eqref{1.0} is
complex valued. Therefore, it is necessary to deal with the multi-valued
nature of the complex logarithm in \eqref{eq1.5}. Even though the 1-D case
was also considered in \cite{KKNS:SIAP2017, KKNS:IP18}, there are three
main differences between the current work and the methods proposed in these publications:

\begin{enumerate}
\item  We propose a simpler way of defining the function $v$ in which the
logarithm is avoided, unlike \cite{KKNS:SIAP2017, KKNS:IP18}. 

\item Item 1 also leads to a coupled system of differential equations of the 
\emph{first order} unlike the ones of the second order in the previous
works. 

\item Item 2 requires, in turn, a different proof of the global strict
convexity of the resulting objective functional.
\end{enumerate}

We refer to \cite{BBE:CPDE2013,BK:NA2015} for similar approaches in the time
domain. The paper \cite{BBE:CPDE2013} is about the reconstruction of the
potential in the wave equation, while \cite{BK:NA2015} is concerned with the
reconstruction of the same coefficient as in the current paper. The
objective functionals in these works are similar to ours, since both of them
use Carleman weight functions, although specific weights are chosen
differently. The main difference between our current work and \cite%
{BBE:CPDE2013,BK:NA2015} is that in those papers at least one initial
condition must be assumed to be nonzero in the entire domain of interest,
whereas we use the delta function as the source. The analysis for the
time-domain problem used in \cite{BBE:CPDE2013, BK:NA2015} cannot be used in
this paper due to the presence of the delta function.

The rest of the paper is organized as follows. In section \ref{sec:2} we
state the forward and inverse problems. Section \ref{sec:3} describes our
version of the method of globally strictly convex functional. The global
strict convexity and the global convergence of the gradient projection
method are discussed in Section \ref{sec:4}. In that section, we also prove
an error estimate for the coefficient to be reconstructed. Section \ref%
{sec:num} discusses some details of the discretization and algorithm.
Numerical results are presented in section \ref{sec:result}.
Finally, concluding remarks are given in Section 7.

\section{Problem statement}

\label{sec:2}

The PDE of the forward problem under consideration is described by equation %
\eqref{1.0}. In this work, we use data at multiple frequencies. Therefore,
in the following we denote the solution of \eqref{1.0} by $u(x,k)$ to
indicate its dependence on the wavenumber. The forward model is then
rewritten as: 
\begin{equation}
u_{xx}(x,k)+k^{2}c(x)u(x,k)=-\delta (x-x^{0}),\ x\in \mathbb{R}.  \label{1.1}
\end{equation}%
In addition, function $u(x,k)$ is assumed to satisfy the following
radiation conditions: 
\begin{equation}
\lim\limits_{x\rightarrow \infty }(u_{x}+iku)=0,\quad
\lim\limits_{x\rightarrow -\infty }(u_{x}-iku)=0.  \label{1.2}
\end{equation}%
Furthermore, we assume that the dielectric constant of the medium is
positive, bounded, and constant outside a given bounded interval $(0,b)$, $%
b>0$. More precisely, the coefficient $c(x)$ is assumed to satisfy: 
\begin{equation}
c\in C^{2}\left( \mathbb{R}\right) ;\ 0<c_{0}\leq c\left( x\right) \leq
1+d,\ \forall x\in \mathbb{R};\ c\left( x\right) =1,\forall x\notin \left(
0,b\right) ,  \label{1.3}
\end{equation}%
where $c_{0}$ and $d$ are given positive numbers. In weak scattering models,
the constant $d$ is usually assumed to be small. However, we do not use this
assumption in this work, i.e., we allow both weak and strong scattering objects. We also assume that the point source $x^{0}$ is
placed outside of the interval where $c\left( x\right) $ is unknown. Without
a loss of generality, we assume throughout of this work that $x^{0}<0.$ The coefficient inverse
problem (CIP) we consider in this paper is stated as follows. \bigskip

\textbf{CIP:} Let $u(x,k)$ be a solution of problem (\ref{1.1})--(\ref{1.2}%
). Suppose that condition (\ref{1.3}) is satisfied. Determine the function $%
c(x)$ for $x\in \left( 0,b\right) ,$ given the following backscatter data 
\begin{equation}
g(k)=u(0,k),\ k\in \lbrack \underline{k},\bar{k}],  \label{1.5}
\end{equation}%
where $[\underline{k},\bar{k}]$ represents the frequency interval used in
the measured data.

\begin{remark}
Since $c(x)=1$ on the interval $(-\infty ,0]$, the Dirichlet data (\ref{1.5}%
) uniquely determines the Neumann data at the same location. Indeed, the
scattered wave $u^{s}:=u-u^{i}$, where $u^{i}$ is the incident wave,
satisfies the Helmholtz equation $u_{xx}^{s}+k^{2}u^{s}=0$ on $(-\infty ,0)$%
, together with the radiation condition $\lim\limits_{x\rightarrow -\infty
}(u_{x}^{s}-iku^{s})=0$. Hence, $u^{s}$ can be written in the form $%
u^{s}(x,k)=Ce^{ikx}.$ The constant $C$ can be calculated from the Dirichlet
data as $C=g(k)-u^{i}(0,k).$ Hence, the Neumann data is given by 
\begin{equation}
g_{1}(k):=u_{x}(0,k)=u_{x}^{i}(0,k)+ikC,\ k\in (\underline{k},\bar{k}).
\label{1.6}
\end{equation}
\end{remark}

\begin{remark}
The uniqueness of this inverse problem has been proved in \cite{KNSN:IPI2016}
under some assumptions about the coefficient $c(x)$. Although these
assumptions are not trivial, we assume in this paper that the uniqueness of
the CIP holds.
\end{remark}

\section{Globally strictly convex functional}

\label{sec:3}

The first idea of this method is to transform problem (\ref{1.1})--(\ref{1.2}%
) into a differential equation which does not contain the
unknown coefficient $c(x)$. After the solution of this equation is found,
the coefficient $c(x)$ can be easily computed. To do that, we define the new
function 
\begin{equation}
v(x,k):=\frac{u_{x}(x,k)}{k^{2}u(x,k)}.  \label{def:v}
\end{equation}%
To guarantee that $v$ is well-defined, we need the following result.

\begin{lemma}
Let $x^0 < 0$ be the position of the point source and $u$ be the solution of problem (\ref{1.1})--(\ref{1.2}).
Under the condition (\ref{1.3}), we have $u\left( x,k\right) \neq 0$ for all 
$x\in \left[ 0,b\right] $ and for all $k>0$.
\end{lemma}

\noindent\textit{Proof}. 
 The proof can be found in \cite{KKNS:IP18}. However, since we need to use some results in this proof in the derivation of the method, we present the proof here. Assume to the contrary that there exists a point $x=a\in \left[ 0,b%
\right] $ and a wavenumber $k_{0}>0$ such that 
\begin{equation}
u\left( a,k_{0}\right) =0.  \label{2.1}
\end{equation}%
Since $c(x) = 1$ for all $x > b$, the solution of  (\ref{1.1})--(\ref%
{1.2}) can be represented as
\begin{equation}
u\left( x,k\right) =C\left( k\right) e^{-ikx},\ \forall x\geq b,\ \forall k>0,
\label{2.2}
\end{equation}%
where $C\left( k\right) $ is a function of $k$. Set in (\ref{1.1}) $k=k_{0},$
multiply this equation by the complex conjugate $\overline{u}\left( 
x,k_{0}\right) $
of $u$ and integrate over the interval $\left( a,b\right) .$ Since $x^0 < 0$, the
right-hand side of the resulting equality is zero. Using (\ref{2.1}), we obtain%
\begin{equation}
\overline{u}\left( b,k_{0}\right) u_{x}\left( b,k_{0}\right)
+\int_{a}^{b}\left[ -\left\vert u_{x}\right\vert
^{2}+k_{0}^{2}c\left( x\right) \left\vert u\right\vert ^{2}\right] dx =0.
\label{2.3}
\end{equation}%
By (\ref{2.2}) $u_{x}\left( b,k_{0}\right) =-ik_{0}u\left( b,k_{0}\right) .$
Hence, $\overline{u}\left( b,k_{0}\right) u_{x}\left( b,k_{0}\right)
=-i k_0\left\vert u\left( b,k_{0}\right) \right\vert ^{2}.$ Hence, (\ref{2.3})
becomes%
\begin{equation}
i\left\vert u\left( b,k_{0}\right) \right\vert ^{2}=\int_{a}^{b} 
\left[ -\left\vert u_{x}\right\vert ^{2} + k_{0}^{2}c\left( x\right) \left\vert
u\right\vert ^{2}\right] dx.  \label{2.4}
\end{equation}%
The left-hand side of (\ref{2.4}) is a purely imaginary number, whereas the
right-hand side is a real number. Therefore, both numbers must be equal to zero.
Hence, $u\left( b,k_{0}\right) =u_{x}\left( b,k_{0}\right) =0.$ By (\ref{1.1}%
)  this means that $u\left( x,k_{0}\right) =0$ for $x\geq
x^{0},$ which is impossible. The proof is complete. $\hfill\square$

\bigskip 

We now derive an equation for $v$. From (\ref{def:v}) we have $u_{x}=k^{2}vu$%
. Differentiating both sides of this identity with respect to $x$, we obtain 
\begin{equation*}
u_{xx} =k^2 (v_x u + v u_x) = k^2 u(v_x + k^2 v^2).
\end{equation*}
Substituting this into (\ref{1.1}), noting that the right-hand side is zero
on the interval $(0,b)$ since $x^0 < 0$, we obtain 
\begin{equation}
v_{x}(x,k)+k^{2}v^{2}(x,k)=-c(x),\ x\in (0,b).  \label{2.5}
\end{equation}%
In addition, function $v(x,k)$ satisfies the following boundary conditions
at $x=0$ and $x=b$: 
\begin{equation}
v(0,k)=v_{0}(k),\quad v(b,k)=-\dfrac{i}{k}.  \label{2.6}
\end{equation}%
Here $v_{0}=\frac{g_{1}(k)}{k^{2}g_{0}(k)} $. The second boundary condition
of \eqref{2.6} is derived from (\ref{2.2}).

If function $v$ is known, then coefficient $c(x)$ can be computed
directly using (\ref{2.5}). However, equation (\ref{2.5}) contains two
unknown functions, $v(x)$ and $c(x)$. Therefore, to find $v$ we eliminate
the unknown coefficient $c(x)$ by taking the derivative of both sides of (%
\ref{2.5}) with respect to $k$. We obtain the following equation: 
\begin{equation}
v_{xk}(x,k)+2kv^{2}(x,k)+2k^{2}v(x,k)v_{k}(x,k)=0,\quad x\in (0,b).
\label{eq:vnew}
\end{equation}

To find function $v$ from \eqref{2.6} and \eqref{eq:vnew}, we use the method
of separation of variables. More precisely, we approximate $v$ via the
following truncated series: 
\begin{equation}
v\left( x,k\right) \approx \sum\limits_{n=1}^{N}v_{n}\left( x\right)
f_{n}\left( k\right) ,  \label{2.10}
\end{equation}%
where $\{f_{n}(k)\}_{k=1}^{\infty }$ is an orthonormal basis in $L_{2}(%
\underline{k},\bar{k})$. Functions $f_{n}(k)$ are real valued and we specify
this basis later. Substituting \eqref{2.10} into \eqref{eq:vnew}, we obtain
the following system: 
\begin{equation}
\sum\limits_{n=1}^{N}v_{n}^{\prime }(x)f_{n}^{\prime
}(k)+\sum\limits_{n=1}^{N}\sum%
\limits_{j=1}^{N}v_{n}(x)v_{j}(x)[2kf_{n}(k)f_{j}(k)+2k^{2}f_{n}(k)f_{j}^{%
\prime }(k)]=0,\quad x\in (0,b).  \label{2.11}
\end{equation}%
To be precise, equation (\ref{2.11}) should be understood as an
approximation of \eqref{eq:vnew} since $v$ is approximated by the truncated sum (\ref{2.10}).
Multiplying both sides of (\ref{2.11}) by $f_{m}(k)$ and integrating over $(%
\underline{k},\bar{k})$, we obtain the following system of coupled quasi-linear
equations for $v_{n}(x)$: 
\begin{equation}
\sum\limits_{n=1}^{N}M_{mn}v_{n}^{\prime
}(x)+\sum\limits_{n=1}^{N}\sum\limits_{j=1}^{N}G_{mnj}v_{n}(x)v_{j}(x)=0,\
m=1,\dots ,N,\ x\in (0,b),  \label{2.12}
\end{equation}%
where the coefficients $M_{mn}$ and $G_{mnj}$ are given by 
\begin{equation}\label{eq:M}
M_{mn}=\int_{\underline{k}}^{\bar{k}}f_{n}^{\prime }(k)f_{m}(k)dk,
\end{equation}%
\begin{equation}\label{eq:G}
G_{mnj}=\int_{\underline{k}}^{\bar{k}%
}[2kf_{n}(k)f_{j}(k)+2k^{2}f_{n}(k)f_{j}^{\prime }(k)]f_{m}(k)dk.
\end{equation}

Using the approximation (\ref{2.10}) for $v(x)$, it follows from (\ref{2.5})  that once functions $v_{n}(x)$, 
$n=1,\dots ,N$, are found, coefficient $c(x)$ is approximated by 
\begin{equation}
c(x)\approx -\sum\limits_{n=1}^{N}v_{n}^{\prime
}(x)f_{n}(k)-k^{2}\sum\limits_{n=1}^{N}\sum%
\limits_{j=1}^{N}v_{n}(x)v_{j}(x)f_{n}(k)f_{j}(k).  \label{eq:c2}
\end{equation}

Note that $v_{n}(x),\ n=1,\dots ,N,$ are complex valued functions. In
numerical implementation, it is more convenient to work with real vectors.
For this purpose, we denote by $V_{n}(x)$ and $V_{n+N}(x)$ the real and
imaginary parts of $v_{n}(x)$ and define the vector-valued real function $%
V(x)$ as $V\left( x\right) =\left( V_{1}(x),\dots ,V_{N}(x),V_{N+1},\dots
,V_{2N}\right) ^{T}$. By separating the real and imaginary parts of %
\eqref{2.12}, we obtain the following equations: 
\begin{eqnarray}
&&\sum\limits_{n=1}^{N}M_{mn}V_{n}^{\prime
}+\sum\limits_{n=1}^{N}\sum%
\limits_{j=1}^{N}G_{mnj}(V_{n}V_{j}-V_{n+N}V_{j+N})=0,
\label{2.122} \\
&&\sum\limits_{n=1}^{N}M_{mn}V_{n+N}^{\prime
}+\sum\limits_{n=1}^{N}\sum%
\limits_{j=1}^{N}G_{mnj}(V_{n}V_{j+N}+V_{n+N}V_{j})=0,
\label{2.123}
\end{eqnarray}%
 for $m=1,\dots ,N$ and $x\in (0,b)$. To be compact, we rewrite these equations in the following vector form 
\begin{equation}\label{2.13}
\tilde{M}V^{\prime }+G\left( V\right) =0,\ x\in (0,b),
\end{equation}%
where $\tilde{M}=%
\begin{bmatrix}
M & 0 \\ 
0 & M%
\end{bmatrix}%
$ is a $2N\times 2N$ block matrix, $M$ is an $N\times N$ matrix with
entries $M_{mn}$ defined by \eqref{eq:M}, and $G(V)=\left( G_{1}(V),\dots
,G_{N}(V),G_{N+1}(V),\dots ,G_{2N}(V)\right) ^{T}$ with 
\begin{eqnarray*}
G_{m}(V)
&=&\sum\limits_{n=1}^{N}\sum%
\limits_{j=1}^{N}G_{mnj}(V_{n}V_{j}-V_{n+N}V_{j+N}), \\
G_{m+N}(V)
&=&\sum\limits_{n=1}^{N}\sum%
\limits_{j=1}^{N}G_{mnj}(V_{n}V_{j+N}+V_{n+N}V_{j})
\end{eqnarray*}%
for $m=1,\dots ,N$. System \eqref{2.13} is coupled with the following boundary conditions: 
\begin{equation}
V(0)=V^{0},\quad V(b)=V^{b},  \label{2.124}
\end{equation}%
where $V^{0}=(V_{1}^{0},\dots ,V_{2N}^{0})^{T}$ and $%
V^{b}=(V_{1}^{b},\dots ,V_{2N}^{b})^{T}$ whose components are
calculated from \eqref{2.6} as follows 
\begin{eqnarray*}
&&V_{n}^{0}=\int_{\underline{k}}^{\bar{k}}\Real(v_{0}(k))f_{n}(k)dk,\quad
V_{n+N}^{0}=\int_{\underline{k}}^{\bar{k}}\Imag(v_{0}(k))f_{n}(k)dk, \\
&&V_{n}^{b}=0,\quad V_{n+N}^{b}=-\int_{\underline{k}}^{\bar{k}}\frac{1}{k}%
f_{n}(k)dk,\quad n=1,\dots ,N.
\end{eqnarray*}

In solving problem (\ref{2.13})--(\ref{2.124}), we require that matrix $M$
be non-singular (and so is $\tilde{M}$). To satisfy this requirement, the
basis $\{f_{n}\}$ must be chosen appropriately. Here we use the same basis of $L_2(\underline{k},\overline{k})$ that was introduced in \cite%
{Kli:JIIP17}. This basis was also used in  \cite{KKNS:IP18}. We start with the complete set $\{k^{n-1}e^{k}\}_{n=1}^{\infty }$
in $L_{2}(0,1)$. Then, we use the Gram-Schmidt orthonormalization procedure
to obtain an orthonormal basis $\{\varphi _{n}(k)\}_{n=1}^{\infty }$ of $%
L_{2}(0,1)$. Finally, we define $f_{n}(k)$ as 
\begin{equation*}
f_{n}(k)=\frac{1}{\sqrt{\bar{k}-\underline{k}}}\varphi _{n}\left( \frac{k-%
\underline{k}}{\bar{k}-\underline{k}}\right) .
\end{equation*}%
It is clear that $\{f_{n}(k)\}_{n=1}^{\infty }$ is an orthonormal basis in $%
L_{2}(\underline{k},\bar{k})$. Moreover, it was proved in \cite{Kli:JIIP17} that
matrix $M$ is upper-triangular with $\det (M)=(\bar{k}-\underline{k})^{-N}$.
Therefore, both matrices $M$ and $\tilde{M}$ are invertible. \bigskip

Next, we introduce a new vector-valued function $Q(x)$ as $Q(x)=V(x)-\hat{V}%
(x)$, where $\hat{V}$ is defined by 
\begin{equation}
\hat{V}(x)=V^{0}+(V^{b}-V^{0})\frac{x}{b},\quad x\in \lbrack 0,b].
\label{Vhat}
\end{equation}%
So each component of $\hat{V}(x)$ is linear on $[0,b]$. The new function $%
Q(x)$ satisfies the following boundary value problem: 
\begin{eqnarray}
&&\left( Q^{\prime }+\hat{V}^{\prime }+F(Q+\hat{V})\right) (x)=0,\ x\in
(0,b),  \label{eq:Q1} \\
&&Q(0)=Q(b)=0.  \label{eq:Q2}
\end{eqnarray}%
where $F(V)(x)=\tilde{M}^{-1}G(V)(x)$. Note that $\hat{V}^{\prime} = \frac{1}{b}(V^{b}-V^{0})$.
Moreover, since $G$ is a quadratic vector function of $V$, so is $F$. 

If the vector function $Q$ is determined, so is $V$. Then, $v(x,k)$ and $c(x)
$ can be calculated using \eqref{2.10} and \eqref{2.5}, respectively.
Therefore, the analysis below focuses on solving the boundary value problem (%
\ref{eq:Q1})--(\ref{eq:Q2}). Let $\mathbb{H}^{1}(0,b)$ be the space of
2N-component vector-valued real functions whose components belong to the
Sobolev space $H^{1}(0,b)$, i.e., $\mathbb{H}^{1}(0,b)=(H^{1}(0,b))^{2N}$.
For each $Q=(Q_{1},Q_{1},\dots ,Q_{2N})^{T}\in \mathbb{H}^{1}(0,b)$, its
norm is defined as 
\begin{equation*}
\Vert Q\Vert _{\mathbb{H}^{1}}:=\left( \sum\limits_{n=1}^{2N}\Vert
Q_{n}\Vert _{H^{1}}^{2}\right) ^{1/2},
\end{equation*}%
where $\Vert \cdot \Vert _{H^{1}}$ denotes the $H^{1}(0,b)$ norm.  We also
define the space $\mathbb{H}_{0}^{1}(0,b):=\{Q(x)\in \mathbb{H}%
^{1}(0,b):Q(0)=Q(b)=0\}$. For each positive number $R$, we denote by $B_{R}$
the ball of radius $R$ centered at the origin in $\mathbb{H}_{0}^{1}(0,b)$,
i.e., 
\begin{equation}
B_{R}:=\{Q\in \mathbb{H}^{1}(0,b):\Vert Q\Vert _{\mathbb{H}^{1}}<R,\
Q(0)=Q(b)=0\}.  \label{G}
\end{equation}

Note that the boundary value problem (\ref{eq:Q1})--(\ref{eq:Q2}) is
over-determined since (\ref{eq:Q1}) is a first order system but there are
two boundary conditions. We also recall that equation \eqref{eq:Q1} is
actually an approximation. Therefore, instead of solving this problem
directly, we approximate $Q$ by minimizing the following Carleman weighted cost functional in the ball $B_R$
\begin{equation}
J_{\lambda ,\alpha }\left( Q\right) =\int_{0}^{b}\left\Vert \left( Q^{\prime
}+\hat{V}^{\prime }+F(Q+\hat{V})\right) (x)\right\Vert _{2N}^{2}e^{-2\lambda
x}dx+\alpha \Vert Q\Vert _{\mathbb{H}^{1}}^{2},  \label{2.16}
\end{equation}%
where $\Vert \cdot \Vert _{2N}$ denotes the Euclidean norm in $\mathbb{R}%
^{2N}$. The exponential term $e^{-2\lambda x}$ is known as the Carleman
Weight Function for the operator $d/dx$.

\begin{remark}
In our theoretical analysis we do not actually need to add the regularization term. However, we keep it here since we have noticed in our numerical studies that we really need it in computations.
\end{remark}

\section{Convexity, global convergence, and accuracy estimate}

\label{sec:4}

In this section, we prove the global strict convexity of the objective
functional $J_{\lambda ,\alpha }$. Next, we prove the global convergence of
the gradient projection method and provide an accuracy estimate of the
reconstructed solution. First, we prove the following Carleman estimate.

\begin{lemma}
\label{le:Carleman} Let $h$ be a real valued function in $H^{1}\left(
0,b\right) $ such that $h\left( 0\right) =0$ and $\lambda $ be a positive
number. Then the following Carleman estimate holds 
\begin{equation}
\int_{0}^{b}\left( h^{\prime }\right) ^{2}e^{-2\lambda x}dx\geq \lambda
^{2}\int_{0}^{b}h^{2}e^{-2\lambda x}dx.
\end{equation}
\end{lemma}

\noindent\textit{Proof}.  Consider the function $v\left( x\right) =h\left( x\right)
e^{-\lambda x}.$ Then $h=v e^{\lambda x}$ and $h^{\prime }=\left( v^{\prime
}+\lambda v\right) e^{\lambda x}.$ Hence,%
\begin{equation*}
\left( h^{\prime }\right) ^{2}e^{-2\lambda x}=\left( v^{\prime }+\lambda
v\right) ^{2}\geq \lambda ^{2}v^{2}+2\lambda v^{\prime }v=\lambda
^{2}v^{2}+\lambda\left(  v^{2}\right) ^{\prime }.
\end{equation*}%
Hence,%
\begin{equation*}
\int_{0}^{b}\left( h^{\prime }\right) ^{2}e^{-2\lambda x}dx\ge
\lambda ^{2}\int_{0}^{b}v^{2}dx+ \lambda \int_{0}^{b}\left( 
v^{2}\right) ^{\prime }dx.
\end{equation*}%
Since 
\begin{equation*}
\lambda \int_{0}^{b}\left( v^{2}\right) ^{\prime }dx= \lambda 
v^{2}(b) -\lambda  v^{2}(0)
=\lambda 
v^{2}(b)\geq 0,
\end{equation*}%
then%
\begin{equation*}
\int_{0}^{b}\left( h^{\prime }\right) ^{2}e^{-2\lambda x}dx\geq
\lambda ^{2}\int_{0}^{b}v^{2}dx=\lambda
^{2}\int_{0}^{b}h^{2}e^{-2\lambda x}dx.
\end{equation*}
The proof is complete. $\hfill\square$
\bigskip 

Next, we prove that the objective functional $J_{\lambda,\alpha}(Q)$ is
smooth.

\begin{lemma}
\label{le2} Let $R$, $\lambda $, and $\alpha $ be arbitrary real numbers
such that $R>0$, $\lambda >0$, and $\alpha \geq 0$. Then, the objective
functional $J_{\lambda ,\alpha }(Q)$ defined by \eqref{2.16} is Fr\'{e}chet
differentiable in $B_{R}$. Moreover, its Fr\'{e}chet derivative $J_{\lambda
,\alpha }^{\prime }(Q)$ is Lipschitz continuous on $B_{R}$, i.e., there
exists a constant $D>0$ depending only on $R$, $F$, $\hat{V}$, $N$, and $%
\alpha $ such that for all $\tilde{Q},\ Q\in B_{R}$ the following inequality
holds 
\begin{equation}
\Vert J_{\lambda ,\alpha }^{\prime }(\tilde{Q})-J_{\lambda ,\alpha }^{\prime
}(Q)\Vert _{\mathbb{H}^{1}}\leq D\Vert \tilde{Q}-Q\Vert _{\mathbb{H}^{1}}.
\label{eq:Lipschitz}
\end{equation}
\end{lemma}

\noindent\textit{Proof}. Since $F$ is a quadratic function, the smoothness of $%
J_{\lambda ,\alpha }$ follows from standard functional analysis arguments.
Indeed, let $\tilde{Q}$ and $Q$ be two functions in $B_{R}$ and denote by $%
h:=\tilde{Q}-Q$. Since $F(Q)$ is a quadratic vector-valued function of $Q$,
it follows that 
\begin{equation}
\begin{split}
\tilde{Q}+\hat{V}^{\prime }+F(\tilde{Q}+\hat{V})& =Q^{\prime }+\hat{V}%
^{\prime }+h^{\prime }+F(Q+h+\hat{V}) \\
& =Q^{\prime }+\hat{V}^{\prime }+F(Q+\hat{V})+L(Q+\hat{V},h)+h^{\prime
}+F(h),
\end{split}
\label{exp}
\end{equation}%
where $L$ is a bilinear operator. Replacing \eqref{exp} into \eqref{2.16},
we have 
\begin{eqnarray}
J_{\lambda ,\alpha }(\tilde{Q}) &=&\int_{0}^{b}e^{-2\lambda x}\Vert \left(
Q^{\prime }+\hat{V}^{\prime }+F(Q+\hat{V})+L(Q+\hat{V},h)+h^{\prime
}+F(h)\right) (x)\Vert _{2N}^{2}dx  \notag \\
&&+\alpha \Vert Q+h\Vert _{\bbH^{1}}^{2}  \notag \\
&=&J_{\lambda ,\alpha }(Q)+\int_{0}^{b}e^{-2\lambda x}\Vert \left( L(Q+\hat{V%
},h)+h^{\prime }+F(h)\right) (x)\Vert _{2N}^{2}dx  \label{eq3.24} \\
&&+2\int_{0}^{b}e^{-2\lambda x}\langle Q^{\prime }+\hat{V}^{\prime }+F(Q+%
\hat{V}),h^{\prime }+(L(Q+\hat{V},h) \rangle_{2N}dx.  \notag \\
&&+2\int_{0}^{b}e^{-2\lambda x}\langle Q^{\prime }+\hat{V}^{\prime }+F(Q+%
\hat{V}),F(h)\rangle_{2N} dx+2\alpha \langle Q,h\rangle_{\bbH^1} +\alpha \Vert h\Vert _{%
\bbH^{1}}^{2}.  \notag
\end{eqnarray}%
Here and below, $\langle ,\rangle_{2N} $ and $\langle, \rangle_{\bbH^1}$ are the inner products in $\mathbb{R}^{2N}$  and in $\bbH^1(0,b)$, respectively. Since $F(Q)$ is a quadratic function of $%
Q$, we have 
\begin{equation*}
\int_{0}^{b}e^{-2\lambda x}\Vert \left( h^{\prime }+L(Q+\hat{V}%
,h)+F(h)\right) (x)\Vert _{2N}^{2}dx=O(\Vert h\Vert _{\bbH^{1}}^{2})
\end{equation*}%
and 
\begin{equation*}
\int_{0}^{b}e^{-2\lambda x}\langle Q^{\prime }+\hat{V}^{\prime }+F(Q+\hat{V}%
),F(h)\rangle_{2N} dx=O(\Vert h\Vert _{\bbH^{1}}^{2})
\end{equation*}%
when $\Vert h\Vert _{\bbH^{1}}\rightarrow 0$. Therefore, it follows from %
\eqref{eq3.24} that 
\begin{equation}
\begin{split}
J_{\lambda ,\alpha }(\tilde{Q})-J_{\lambda ,\alpha }(Q)&
=2\int_{0}^{b}e^{-2\lambda x}\langle Q^{\prime }+\hat{V}^{\prime }+F(Q+\hat{V%
}),h^{\prime }+L(Q+\hat{V},h)\rangle_{2N} dx \\
& +2\alpha \langle Q,h\rangle_{\bbH^1} +O(\Vert h\Vert _{\bbH^{1}}^{2}).
\end{split}
\label{eq3.25}
\end{equation}%
Since the first two terms on the right-hand side of \eqref{eq3.25} are
bounded linear operators of $h$, we conclude that $J_{\lambda ,\alpha }$ is
Fr\'{e}chet differentiable and its gradient is given by 
\begin{equation}
J_{\lambda ,\alpha }^{\prime }(Q)\left( h\right) =2\int_{0}^{b}e^{-2\lambda
x}\langle Q^{\prime }+\hat{V}^{\prime }+F(Q+\hat{V}),h^{\prime }+L(Q+\hat{V}%
,h)\rangle_{2N} dx+2\alpha \langle Q,h\rangle_{\bbH^1}
\label{gradJ}
\end{equation}%
for any $h\in \bbH^{1}(0,b)$.

Next, we prove the Lipschitz continuity of $J_{\lambda ,\alpha }^{\prime }$.
For an arbitrary vector-valued function $h\in \bbH^{1}(0,b)$, \eqref{gradJ}
implies that 

\begin{equation}
\begin{split}
& \left[ J_{\lambda ,\alpha }^{\prime }(\tilde{Q})-J_{\lambda ,\alpha
}^{\prime }(\tilde{Q})\right] \left( h\right)  \\
& =2\int_{0}^{b}e^{-2\lambda x}\langle \tilde{Q}^{\prime }+\hat{V}^{\prime
}+F(\tilde{Q}+\hat{V}),h^{\prime }+L(\tilde{Q}+\hat{V},h)\rangle_{2N} dx \\
& -2\int_{0}^{b}e^{-2\lambda x}\langle Q^{\prime }+\hat{V}^{\prime }+F(Q+%
\hat{V}),h^{\prime }+L(Q+\hat{V},h)\rangle_{2N} dx \\
& +2\alpha \langle Q,h \rangle_{\bbH^1} .
\end{split}
\label{gradJ2}
\end{equation}%
To estimate the first two terms on the right-hand side of \eqref{gradJ2}, we separate them as
follows: 
\begin{equation}
\begin{split}
2& \int_{0}^{b}e^{-2\lambda x}\langle \tilde{Q}^{\prime }+\hat{V}^{\prime
}+F(\tilde{Q}+\hat{V}),h^{\prime }+L(\tilde{Q}+\hat{V},h)\rangle_{2N} dx \\
-& 2\int_{0}^{b}e^{-2\lambda x}\langle Q^{\prime }+\hat{V}^{\prime }+F(Q+%
\hat{V}),h^{\prime }+L(Q+\hat{V},h)\rangle_{2N} dx \\
=& 2\int_{0}^{b}e^{-2\lambda x}\langle \tilde{Q}^{\prime }-Q^{\prime }+F(%
\tilde{Q}+\hat{V})-F(Q+\hat{V}),h^{\prime }+L(\tilde{Q}+\hat{V},h)\rangle_{2N} dx
\\
+& 2\int_{0}^{b}e^{-2\lambda x}\langle Q^{\prime }+\hat{V}^{\prime }+F(Q+%
\hat{V}),L(\tilde{Q}-Q,h)\rangle_{2N} dx.
\end{split}
\label{gradJ22}
\end{equation}%
In obtaining the last term, we have used the identity $L(\tilde{Q}+\hat{V}%
,h)-L(Q+\hat{V},h)=L(\tilde{Q}-Q,h)$ thanks to the bilinearity of $L$. Since $F$ is a quadratic function of $Q
$, there exist constants $D_{1}=D_{1}(R,F,\hat{V},N)$ and $D_{2}=D_{2}(R,F,%
\hat{V},N)$ such that for all $\tilde{Q},Q\in B_{R}$ and for all $x\in (0,b)$, 
\begin{equation}
\begin{split}
& \Vert (\tilde{Q}^{\prime }-Q^{\prime }+F(\tilde{Q}+\hat{V})-F(Q+\hat{V}))(x)\Vert
_{2N}^{2}\\
& =\Vert (\tilde{Q}^{\prime }-Q^{\prime }+F(\tilde{Q}-Q)+L(\tilde{Q}%
-Q,Q+\hat{V}))(x)\Vert _{2N}^{2} \\
& \leq D_{1}^{2}\left( \Vert (\tilde{Q}^{\prime }-Q')(x)\Vert_{2N}^2+\Vert 
(\tilde{Q}-Q)(x)\Vert_{2N}^{2}\right) 
\end{split}
\label{eq4.19}
\end{equation}%
and 
\begin{equation}
\Vert (Q^{\prime }+\hat{V}^{\prime }+F(Q+\hat{V}))(x)\Vert _{2N}\leq D_{2}.
\label{eq4.20}
\end{equation}

On the other hand, there exist constants $D_3 = D_3(R,F,\hat V, N)$ and $D_4
= D_4(R,F,\hat V, N)$ such that  
\begin{eqnarray}
&& \| (h^{\prime }+ L(\tilde Q+\hat V, h))(x)\|^2_{2N} \le D_3^2
\left(\|h^{\prime}(x)\|^2_{2N} + \|h(x)\|^2_{2N}\right),  \label{eq4.21} \\
&& \|( L(\tilde Q-Q, h))(x)\|_{2N} \le D_4 \|(\tilde Q - Q)(x) \|_{2N} \|h(x)\|_{2N}
\label{eq4.22}
\end{eqnarray}

Hence, using Cauchy-Schwarz inequality, \eqref{eq4.19} and \eqref{eq4.21},
we obtain 
\begin{eqnarray}
&& \Big| \int_0^b e^{-2\lambda x} \langle \tilde Q^{\prime }- Q^{\prime }+
F(\tilde Q + \hat V) - F( Q + \hat V), h^{\prime }+ L(\tilde Q+\hat V, h)
\rangle_{2N} dx \Big|  \notag \\
&& \le \int_0^b \| (\tilde Q^{\prime }- Q^{\prime }+ F(\tilde Q + \hat V) -
F( Q + \hat V))(x)\|_{2N} \|(h^{\prime }+ L(\tilde Q+\hat V, h))(x) \|_{2N} dx  \notag
\\
&& \le \left(\int_0^b \| (\tilde Q^{\prime }- Q^{\prime }+ F(\tilde Q + \hat
V) - F( Q + \hat V))(x)\|^2_{2N} dx\right)^{1/2}  \notag \\
&& \times \left( \int_0^b \|(h^{\prime }+ L(\tilde Q+\hat V, h))(x) \|^2_{2N}
dx\right)^{1/2}  \notag \\
&& \le D_1 D_3 \|\tilde Q - Q\|_{\bbH^1} \|h\|_{\bbH^1}.  \label{eq4.23}
\end{eqnarray}
Note that we have used the fact that $e^{-2\lambda x}\in \left( 0,1\right)$ in the above estimate. Similarly, we have from %
\eqref{eq4.20} and \eqref{eq4.22}  that
\begin{eqnarray}
&&\Big|\int_{0}^{b}e^{-2\lambda x}\langle Q^{\prime }+\hat{V}^{\prime }+F(Q+\hat{%
V}),L(\tilde{Q}+\hat{V},h)-L(Q+\hat{V},h)\rangle_{2N} dx \Big|  \notag \\
&\leq &\int_{0}^{b}\Vert (Q^{\prime }+\hat{V}^{\prime }+F(Q+\hat{V}))(x)\Vert
_{2N}\Vert (L(\tilde{Q}-Q,h))(x)\Vert _{2N}dx  \notag \\
&\leq &\int_{0}^{b}D_{2}D_{4}\Vert (\tilde{Q}-Q)(x)\Vert _{2N}\Vert h(x)\Vert
_{2N}dx\leq D_{2}D_{4}\Vert \tilde{Q}-Q\Vert _{\bbH^{1}}\Vert h\Vert _{\bbH%
^{1}}.  \label{eq4.24}
\end{eqnarray}%
It follows from \eqref{gradJ2}, \eqref{gradJ22}, \eqref{eq4.23}, and %
\eqref{eq4.24} that 
\begin{equation*}
\left\vert \left[ J_{\lambda ,\alpha }^{\prime }(\tilde{Q})-J_{\lambda
,\alpha }^{\prime }(Q)\right] \left( h\right) \right\vert \leq
2(D_{1}D_{3}+D_{2}D_{4}+\alpha )\Vert \tilde{Q}-Q\Vert _{\bbH^{1}}\Vert
h\Vert _{\bbH^{1}}.
\end{equation*}%
This inequality implies \eqref{eq:Lipschitz}. The proof is complete. $\hfill \square$ 
\bigskip 

We are now ready to state and prove our main theoretical results.

\begin{theorem}[\textbf{Convexity}]
\label{th1} Let $R$ be an arbitrary positive number and $\alpha \geq 0$.
Then, there exists a sufficiently large number $\lambda _{0}=\lambda
_{0}( R,F,\hat{V},N) >0$ such that the objective functional $%
J_{\lambda ,\alpha }(Q)$ defined by \eqref{2.16} is strictly convex on $%
\overline{B}_{R}$ for all $\lambda \geq \lambda _{0}$. More precisely, there
exists a constant $C^{\ast }=C^{\ast }( R,F,\hat{V},N) >0$ such
that for arbitrary vector functions $\tilde{Q},\ Q\in \overline{B}_{R}$, the
following estimate holds: 
\begin{equation}
J_{\lambda ,\alpha }(\tilde{Q})-J_{\lambda ,\alpha }(Q)-J_{\lambda ,\alpha
}^{\prime }(Q)(\tilde{Q}-Q)\geq C^{\ast }\Vert \tilde{Q}-Q\Vert _{\mathbb{H}%
^{1}}^{2},\text{ }\forall \lambda \geq \lambda _{0}.  \label{convex}
\end{equation}%
Both constants $\lambda _{0}$ and $C^{\ast }$ depend only on the listed
parameters.
\end{theorem}

\noindent\textit{Proof.} Denote $h=\tilde{Q}-Q$. It follows from Lemma \ref{le2} that $%
J_{\lambda ,\alpha }(Q)$ is Fr\'{e}chet differentiable on $B_{2R}$. It
follows from \eqref{eq3.24} and \eqref{gradJ} that 
\begin{equation}
\begin{split}
& J_{\lambda ,\alpha }(\tilde{Q})-J_{\lambda ,\alpha }(Q)-J_{\lambda ,\alpha
}^{\prime }(Q)h \\
& =\int_{0}^{b}e^{-2\lambda x}\Vert ( L(Q+\hat{V},h)+h^{\prime
}+F(h)) (x)\Vert _{2N}^{2}dx \\
& +2\int_{0}^{b}e^{-2\lambda x}\langle Q^{\prime }+\hat{V}^{\prime }+F(Q+%
\hat{V}),F(h)\rangle_{2N} dx+\alpha \Vert h\Vert _{\bbH^{1}}^{2}.
\end{split}
\label{eq3.29}
\end{equation}%
Here $L$ is the same bilinear operator as in Lemma \ref{le2}. The first term
on the right-hand side of \eqref{eq3.29} is estimated as follows: 
\begin{eqnarray}
&&\Vert (L(Q+\hat{V},h)+h^{\prime }+F(h))(x)\Vert _{2N}^{2} \notag \\
&=&\Vert h^{\prime }(x)\Vert _{2N}^{2}+2\langle h^{\prime }(x),(L(Q+\hat{V}%
,h)+F(h))(x)\rangle_{2N} +\Vert (L(Q+\hat{V},h)+F(h))(x)\Vert _{2N}^{2}  \notag \\
&\geq &\Vert h^{\prime }(x)\Vert _{2N}^{2}-\frac{1}{2}\Vert h^{\prime }(x)\Vert
_{2N}^{2}-2\Vert (L(Q+\hat{V},h)+F(h))(x)\Vert _{2N}^{2} \notag\\
& + & \Vert (L(Q+\hat{V}%
,h)+F(h))(x)\Vert _{2N}^{2}  \notag \\
&=&\frac{1}{2}\Vert h^{\prime }(x)\Vert _{2N}^{2}-\Vert (L(Q+\hat{V}%
,h)+F(h))(x)\Vert _{2N}^{2}.  \label{eq3.30}
\end{eqnarray}%
In addition, since $F(h)$ is a quadratic vector function of $h$ and $L$ is a
bilinear operator, there is a constant $C_{1}=C_{1}(R,F,\hat{V},N)$ such
that 
\begin{equation*}
\Vert (L(Q+\hat{V},h)+F(h))(x)\Vert _{2N}^{2}\leq C_{1}\Vert h(x)\Vert _{2N}^{2}.
\end{equation*}%
Substituting this inequality into \eqref{eq3.30}, we obtain 
\begin{equation}
\Vert (L(Q+\hat{V},h)+h^{\prime }+F(h))(x)\Vert _{2N}^{2}\geq \frac{1}{2}\Vert
h^{\prime }(x)\Vert _{2N}^{2}-C_{1}\Vert h(x)\Vert _{2N}^{2}.  \label{eq4.13}
\end{equation}

On the other hand, since $\Vert Q\Vert _{\bbH^{1}}\leq R$, the second term
on the right-hand side of \eqref{eq3.29} is estimated as: 
\begin{equation}
2\langle (Q^{\prime }+\hat{V}^{\prime }+F(Q+\hat{V}))(x),F(h)(x)\rangle_{2N} \geq
-C_{2}\Vert h(x)\Vert _{2N}^{2},  \label{eq3.31}
\end{equation}%
where $C_{2}=C_{2}(R,F,\hat{V},N)$ is a constant depending only on $R$, $F$, 
$\hat{V}$, and $N$. Substituting \eqref{eq4.13} and \eqref{eq3.31} into %
\eqref{eq3.29}, we obtain 
\begin{equation}
\begin{split}
& J_{\lambda ,\alpha }(\tilde{Q})-J_{\lambda ,\alpha }(Q)-J_{\lambda ,\alpha
}^{\prime }(Q)h \\
& \geq \frac{1}{2}\int_{0}^{b}e^{-2\lambda x}\Vert h^{\prime }(x)\Vert
_{2N}^{2}dx-C_{3}\int_{0}^{b}e^{-2\lambda x}\Vert h(x)\Vert _{2N}^{2}dx+\alpha
\Vert h\Vert _{\bbH^{1}}^{2},
\end{split}
\label{eq3.32}
\end{equation}%
where $C_{3}=C_{1}+C_{2}$. Let $\lambda _{0}$ be a positive constant such
that $\lambda _{0}^{2}/8>C_{3}$. Note that $h\in \bbH^{1}(0,b)$ and $h(0)=0$%
. Applying Lemma \ref{le:Carleman} to the first term in the second line of %
\eqref{eq3.32}, we obtain
\begin{equation}
\begin{split}
& J_{\lambda ,\alpha }(\tilde{Q})-J_{\lambda ,\alpha }(Q)-J_{\lambda ,\alpha
}^{\prime }(Q)\left( h\right)  \\
& \geq \frac{1}{4}\int_{0}^{b}e^{-2\lambda x}\Vert h^{\prime }\Vert
_{2N}^{2}dx+\frac{\lambda _{0}^{2}}{8}\int_{0}^{b}e^{-2\lambda x}\Vert
h\Vert _{2N}^{2}dx+\alpha \Vert h\Vert _{\bbH^{1}}^{2} \\
& \geq (C_{4}+\alpha )\Vert h\Vert _{\bbH^{1}}^{2},
\end{split}
\label{eq4.15}
\end{equation}%
for all $\lambda \geq \lambda _{0}$. Here $C_{4}=e^{-2\lambda b}\min
\{1/4,\lambda _{0}^{2}/8\}$. Setting $C^{\ast }=e^{-2\lambda _{0}b}\min
\{1/4,\lambda _{0}^{2}/8\} + \alpha$, we obtain \eqref{convex}. The proof is complete. $\hfill\square$ 
\bigskip 

%
%
%

Due to the convexity and smoothness of $J_{\lambda ,\alpha }$ on the closed
ball $\overline{B}_{R}$ the following result follows from theorem 2.1 of 
\cite{BKK}:

\begin{theorem}
\label{th3} Assume that the parameter $\lambda \geq \lambda _{0}$ is chosen
according to Theorem \ref{th1}. Then, the objective functional $J_{\lambda
,\alpha }$ has a unique minimizer $Q^{(min)}$ on $\overline{B}_{R}.$
Furthermore, the following condition holds true: 
\begin{equation*}
J_{\lambda ,\alpha }^{\prime }(Q^{(min)})(Q^{(min)}-Q)\leq 0\quad \text{for
all }\ Q\in \overline{B}_{R}.
\end{equation*}
\end{theorem}

To find the minimizer of $J_{\lambda ,\alpha }$ on $\overline{B}_{R}$, we
use gradient-based approaches. One simple method is the gradient
projection algorithm which starts from an initial guess $Q^{(0)}$ in $B_{R}$
and finds approximations of $Q^{(min)}$ using the following iterative
process: 
\begin{equation}
Q^{(n+1)}=\mathcal{P}\left[ Q^{(n)}-\gamma J_{\lambda ,\alpha }^{\prime
}(Q^{(n)})\right] ,\quad n=0,1,\dots   \label{gpa}
\end{equation}%
Here $\mathcal{P}:\mathbb{H}^{1}\rightarrow \overline{B}_{R}$ is the
orthogonal projection operator from $\mathbb{H}^{1}$ onto $\overline{B}_{R}$
and $\gamma $ is the step length at each iteration.

The following theorem ensures that the gradient projection algorithm is
convergent for an arbitrary initial guess $Q^{(0)}\in B_{R}$.

\begin{theorem}[\textbf{Global convergence}]
\label{th4} Assume that the parameter $\lambda \geq \lambda _{0}$ is chosen
according to Theorem \ref{th1}. Let $Q^{(min)}$ be the unique global minimum
of $J_{\lambda ,\alpha }(Q)$ in the closed ball $\overline{B}_{R}$. Let $%
\{Q^{(n)}\}_{n=0}^{\infty }$ be the sequence (\ref{gpa}) of the gradient
projection algorithm \eqref{gpa} with $\gamma $ small enough. Then there
exists a sufficiently small number $\gamma _{0}=\gamma _{0} (R,F,\hat{V}%
,N) \in \left( 0,1\right) $ depending only on the listed parameters such
that for every $\gamma \in \left( 0,\gamma _{0}\right) $ the sequence $%
\{Q^{(n)}\}_{n=0}^{\infty }$ converges to $Q^{(min)}$ in the space $\mathbb{H%
}^{1}.$ Furthermore, there exists a number $q=q\left( \gamma \right) \in
(0,1)$ such that 
\begin{equation}
\Vert Q^{(n)}-Q^{(min)}\Vert _{\mathbb{H}^{1}}\leq q^{n}\Vert
Q^{(0)}-Q^{(min)}\Vert _{\mathbb{H}^{1}},\text{ }n=1,2,...  \label{error1}
\end{equation}
\end{theorem}

\noindent\textit{Proof}. It is sufficient to show that the operator in the right hand
side of (\ref{gpa}) is contractual mapping on $\overline{B}_{R}.$ Note that $%
\mathcal{P}$ maps $\overline{B}_{R}$ into itself. Let $Y,Z\in \overline{B}%
_{R}$ be two arbitrary points of the closed ball $\overline{B}_{R}.$ We have%
\begin{eqnarray}
&&\left\Vert \mathcal{P}\left[ Y-\gamma J_{\lambda ,\alpha }^{\prime }(Y)%
\right] -\mathcal{P}\left[ Z-\gamma J_{\lambda ,\alpha }^{\prime }(Z)\right]
\right\Vert _{\mathbb{H}^{1}}^{2} \notag \\
&&=\left\Vert \left( Y-Z\right) -\gamma \mathcal{P}\left[ J_{\lambda ,\alpha
}^{\prime }(Y)-J_{\lambda ,\alpha }^{\prime }(Z)\right] \right\Vert _{%
\mathbb{H}^{1}}^{2}  \notag\\
&&=\left\Vert Y-Z\right\Vert _{\mathbb{H}^{1}}^{2}+\gamma ^{2}\left\Vert 
\mathcal{P}\left[ J_{\lambda ,\alpha }^{\prime }(Y)-J_{\lambda ,\alpha
}^{\prime }(Z)\right] \right\Vert _{\mathbb{H}^{1}}^{2} \notag\\ 
&&-2\gamma \langle
J_{\lambda ,\alpha }^{\prime }(Y)-J_{\lambda ,\alpha }^{\prime
}(Z),Y-Z\rangle_{\bbH^1} \label{3}.
\end{eqnarray}%
By Lemma 4.2%
\begin{equation}
\gamma ^{2}\left\Vert \mathcal{P}\left[ J_{\lambda ,\alpha }^{\prime
}(Y)-J_{\lambda ,\alpha }^{\prime }(Z)\right] \right\Vert _{\mathbb{H}%
^{1}}^{2}=\gamma ^{2}\left\Vert J_{\lambda ,\alpha }^{\prime }(Y)-J_{\lambda
,\alpha }^{\prime }(Z)\right\Vert _{\mathbb{H}^{1}}^{2}\leq \gamma
^{2}D^{2}\left\Vert Y-Z\right\Vert _{\mathbb{H}^{1}}^{2}.  \label{4}
\end{equation}%
We now rewrite (\ref{convex}) in two different forms:%
\begin{equation*}
J_{\lambda ,\alpha }\left( Y\right) -J_{\lambda ,\alpha }\left( Z\right)
-J_{\lambda ,\alpha }^{\prime }(Z)\left( Y-Z\right) \geq C^{\ast }\left\Vert
Y-Z\right\Vert _{\mathbb{H}^{1}}^{2},
\end{equation*}%
\begin{equation*}
J_{\lambda ,\alpha }\left( Z\right) -J_{\lambda ,\alpha }\left( Y\right)
+J_{\lambda ,\alpha }^{\prime }(Y)\left( Y-Z\right) \geq C^{\ast }\left\Vert
Y-Z\right\Vert _{\mathbb{H}^{1}}^{2}.
\end{equation*}%
Summing up, we obtain%
\begin{equation*}
\langle J_{\lambda ,\alpha }^{\prime }(Y)-J_{\lambda ,\alpha }^{\prime
}(Z),Y-Z\rangle_{\bbH^1} \geq 2C^{\ast }\left\Vert Y-Z\right\Vert _{\mathbb{H}%
^{1}}^{2}.
\end{equation*}%
Substituting this into (\ref{3}) and using (\ref{4}), we obtain%
\begin{equation*}
\left\Vert \mathcal{P}\left[ Y-\gamma J_{\lambda ,\alpha }^{\prime }(Y)%
\right] -\mathcal{P}\left[ Z-\gamma J_{\lambda ,\alpha }^{\prime }(Z)\right]
\right\Vert _{\mathbb{H}^{1}}^{2}\leq \left( 1-4C^{\ast }\gamma +\gamma
^{2}D^{2}\right) \left\Vert Y-Z\right\Vert _{\mathbb{H}^{1}}^{2}.
\end{equation*}%
If $0<\gamma <4C^{\ast }/D^{2},$ then the number $q\left( \gamma \right)
=\left( 1-4C^{\ast }\gamma +\gamma ^{2}D^{2}\right) \in \left( 0,1\right) .$
The proof is complete. $\hfill\square$
\bigskip 

Finally, we discuss the reconstruction accuracy. For this purpose, denote by 
$\hat{V}^{\ast }(k)$ the function defined by \eqref{Vhat} with $V^{0}$
replaced by the noise-free data $V^{0,\ast }$ at $x=0$. Assume that $\Vert
V^{0,\ast }-V^{0}\Vert _{2N}<\delta $. It follows from \eqref{Vhat} that
there exists a constant $\beta >0$ depending only on $b$ such that 
\begin{equation}
\Vert \hat{V}-\hat{V}^{\ast }\Vert _{\mathbb{H}^{1}}<\beta \delta .
\label{noise}
\end{equation}

%

We have the following result for error estimates.

\begin{theorem}[\textbf{Error estimates}]
\label{th5} Assume that there exists an exact solution $Q^{\ast }$ of
problem \eqref{eq:Q1}--\eqref{eq:Q2} in $\overline{B}_{R}$ associated with
the exact data $V^{0,\ast }$ and the coefficient $c^{\ast }(x)$ is
calculated from \eqref{eq:c2} with $V=V^{\ast }:=Q^{\ast }+\hat{V}^{\ast }$.
Let $Q^{(min)}$ be the unique minimizer of the objective functional $%
J_{\lambda ,\alpha }$ on $\overline{B}_{R}$ and $c^{(min)}$ be calculated
from \eqref{eq:c2} with $V=Q^{(min)}+\hat{V}$. Let $\lambda _{0}$ be chosen
as in Theorem \ref{th1}. Moreover, $\alpha = \xi\delta ^{2}$ for some constant $\xi > 0$. Then, for $\lambda
\geq \lambda _{0}$, we have the following error estimates: 
\begin{eqnarray}
&&\Vert Q^{\ast }-Q^{(min)}\Vert _{\mathbb{H}^{1}}\leq C^{\ast \ast }\delta ,
\label{errorQ} \\
&&\Vert c^{\ast }-c^{(min)}\Vert _{L^{2}}\leq C^{\ast \ast }\delta ,
\label{errorC}
\end{eqnarray}%
where the number $C^{\ast \ast }=C^{\ast \ast } ( F,R,N) >0$
depends only on the listed parameters.
\end{theorem}

\noindent\textit{Proof}. Since $Q^{\ast }$ and $Q^{(min)}$ belong to $\overline{B}_{R}$%
, Theorems \ref{th1} and \ref{th3} imply that 
\begin{eqnarray}
J_{\lambda ,\alpha }(Q^{\ast }) &\geq &J_{\lambda ,\alpha }(Q^{\ast
})-J_{\lambda ,\alpha }(Q^{(min)})-J_{\lambda ,\alpha }^{\prime
}(Q^{(min)})(Q^{\ast }-Q^{(min)})  \notag \\
&\geq &C^{\ast }\Vert Q^{\ast }-Q^{(min)}\Vert _{\bbH^{1}}^2.  \label{err1}
\end{eqnarray}%
On the other hand, it follows from \eqref{eq:Q1} that 
\begin{eqnarray*}
&&(Q^{\ast})' + \hat{V}'+F( Q^{\ast}+\hat{V}) \\
&=&\left[ (Q^{\ast})'+(\hat{V}^{\ast})'+F(Q^{\ast }+\hat{V}^{\ast
})\right] +\left[ \hat{V}^{\prime }-(\hat{V}^{\ast})'+F(Q^{\ast }+\hat{V%
})-F(Q^{\ast }+\hat{V}^{\ast })\right]  \\
&=&\hat{V}^{\prime }-(\hat{V}^{\ast})'+F(Q^{\ast }+\hat{V%
})-F(Q^{\ast }+\hat{V}^{\ast }).
\end{eqnarray*}%
Obviously 
\begin{equation*}
\int_{0}^{b}\left\Vert \left( \hat{V}^{\prime }-\hat{V}^{\ast \prime
}+F(Q^{\ast }+\hat{V})-F(Q^{\ast }+\hat{V}^{\ast })\right) \left( x\right)
\right\Vert _{2N}^{2}dx\leq C_{5}\delta ^{2},
\end{equation*}%
where the constant $C_{5}=C_{5}(F,R,N) >0$ depends only on
the listed parameters. 
Hence, 
\begin{eqnarray}
J_{\lambda ,\alpha }(Q^{\ast }) &=&\int_{0}^{b}\left\Vert \left( \hat{V}%
^{\prime }-\hat{V}^{\ast \prime }+F(Q^{\ast }+\hat{V})-F(Q^{\ast }+\hat{V}%
^{\ast })\right) \left( x\right) \right\Vert _{2N}^{2}e^{-2\lambda x}dx 
\notag \\
&+&\alpha \Vert Q^{\ast }\Vert _{\bbH^{1}}^{2}\leq C_5
\delta ^{2} +\alpha R^{2} = ( C_5 + \xi R^2) \delta^2.  \label{err2}
\end{eqnarray}%
The error estimate \eqref{errorQ} follows
from \eqref{err1} and \eqref{err2} with $C^{\ast \ast }=\sqrt{(C_5 + \xi R^2)/C^*}$. The error estimate \eqref{errorC} for
coefficient $c(x)$ follows directly from \eqref{errorQ} and \eqref{eq:c2}.
The proof is complete. $\hfill\square$

\begin{remark}
Error estimates for $Q^{(n)}$ and $c^{(n)}$easily follows from Theorems \ref%
{th4}--\ref{th5}.
\end{remark}

\begin{remark}\label{remark2}
Due to the fact that \eqref{eq:c2} is only an approximation, the
\textquotedblleft exact\textquotedblright\ coefficient $c^{\ast }$ in
Theorem \ref{th5} is actually not the true coefficient of the original
inverse problem. The difference between this coefficient and the true one
depends on the truncation error in \eqref{2.10}, which is hard to estimate
analytically. In our numerical analysis presented in Section \ref{sec:result}%
, we demonstrate numerically that this error is small even when only a few
Fourier coefficients are used in \eqref{2.10}.
\end{remark}

\section{Discretization and algorithm}

\label{sec:num} In this section, we describe the discretization and
numerical algorithm for finding the vector function $Q$. For the numerical
implementation, it is more convenient to use \eqref{2.13} than \eqref{eq:Q1} because all coefficients in \eqref{2.13} are explicitly given. Note that the two forms are equivalent. In addition, suppose that the measured data $g(k)=u(0,k)$ is available at a finite number of 
wavenumbers $k=k_{1},\dots ,k_{K}$. In this case, we consider each basis function $%
\{f_{n}\}_{n=1}^{N}$ as a $K$-dimensional vectors instead of a function in $L_2(\underline{k},\overline{k})$ and replace the $L_{2}(%
\underline{k},\overline{k})$ norm by the inner
product of real valued $K$-dimensional vectors.

\subsection{Discretization with respect to $x$}

We consider a partition of the interval $(0,b)$ into $M$ sub-intervals by a
uniform grid $0=x_{0}<x_{1}<\cdots <x_{M}=b$ with $x_{m+1}-x_{m}=h$, $%
m=0,\dots ,M-1$. We define the discrete variables as $
Q_{h}:=\{Q_{n,m},\ n=1,\dots ,2N;\ m=0\dots ,M\}$ with $Q_{n,m}=Q_{n}(x_{m})$%
. Note that $Q_{n,0}=Q_{n,M}=0$ due to \eqref{eq:Q2}. The discrete
approximation $\hat{V}_{h}$ of $\hat{V}$ is defined in the same way. We also
define $V_{h}=Q_{h}+\hat{V}_{h}$. The regularized discrete objective function is written
as: 
\begin{equation}
J_{h,\lambda ,\alpha }(Q_{h}):=h\sum\limits_{n=1}^{N}\sum\limits_{m=0}^{M-1}
\left[ (J_{n,m}^{(1)}(V_{h}))^{2}+(J_{n,m}^{(2)}(V_{h}))^{2}\right] \varphi
_{m} + \alpha \mathcal{R}(Q_h),  \label{eq:numqn}
\end{equation}%
where $\varphi _{m}=e^{-2\lambda x_{m}}$, the function $\mathcal R(Q_h)$ is the regularization term given by
\begin{equation}\label{eq:reg}
\mathcal{R}(Q_h) = h\sum\limits_{n=1}^{2N}\sum\limits_{m=0}^{M-1} \left[Q_{n,m}^{2}  + \left(\frac{Q_{n,m+1} - Q_{n,m}}{h}\right)^2\right],
\end{equation}
and the functions $J_{n,m}^{(1)}(Q_{h})$ and $J_{n,m}^{(2)}(Q_{h})$ are defined by
\begin{equation}
J_{n,m}^{(1)}(Q_{h})=\sum\limits_{l=1}^{N}M_{nl}\frac{V_{n,m+1}-V_{n,m}}{h}%
+\sum\limits_{l=1}^{N}\sum%
\limits_{j=1}^{N}G_{nlj}[V_{l,m}V_{j,m}-V_{l+N,m}V_{j+N,m}],
\label{eq:numqn2}
\end{equation}%
\begin{equation}
J_{n,m}^{(2)}(Q_{h})=\sum\limits_{l=1}^{N}M_{nl}\frac{V_{n+N,m+1}-V_{n+N,m}}{%
h}+\sum\limits_{l=1}^{N}\sum%
\limits_{j=1}^{N}G_{nlj}[V_{l,m}V_{j+N,m}+V_{l+N,m}V_{j,m}],
\label{eq:numqn22}
\end{equation}%
for $n=1,\dots ,N$ and $m=0,\dots ,M-1$.
\bigskip 

The unknown variables to be found are $Q_{n,m}, \ n = 1\dots, 2N, \ m =
1,\dots, M-1$. The gradient of the discrete cost function $%
J_{h,\lambda,\alpha}$ can be derived from (\ref{eq:numqn})--(\ref{eq:numqn22}%
). More precisely, using direct calculations, we obtain 
\begin{equation}  \label{eq:grad1}
\frac{\partial J_{h,\lambda,\alpha}(Q_h)}{\partial Q_{s,t}} =
2h\sum\limits_{n=1}^{N}\sum\limits_{m=0}^{M-1}\left[ J^{(1)}_{n,m}(Q_{h}) 
\frac{\partial J^{(1)}_{n,m}(V_h)}{\partial Q_{s,t}} + J^{(2)}_{n,m}(Q_{h}) 
\frac{\partial J^{(2)}_{n,m}(V_h)}{\partial Q_{s,t}} \right] \varphi_m + \alpha \frac{\partial \mathcal{R}(Q_h)}{\partial Q_{s,t}},
\end{equation}%
for $s = 1,\dots, 2N; \ t = 1,\dots, M-1$, where 
\begin{equation*}  \label{eq:grad3}
\frac{\partial J^{(1)}_{n,m}(V_h)}{\partial Q_{s,t}} = \frac{\partial
J^{(2)}_{n,m}(V_h)}{\partial Q_{s+N,t}}= 
\begin{cases}
-\frac{M_{ns}}{h} + \sum\limits_{l = 1}^{N} (G_{nls} + G_{nsl} ) V_{l,m}, & 
t = m,\ s = 1,\dots, N, \\ 
\frac{M_{ns}}{h}, & t = m+1,\ s= 1,\dots, N, \\ 
0, & \text{ otherwise},%
\end{cases}%
\end{equation*}
and 
\begin{equation*}  \label{eq:grad5}
\frac{\partial J^{(2)}_{n,m}(V_h)}{\partial Q_{s,t}} = -\frac{\partial
J^{(1)}_{n,m}(V_h)}{\partial Q_{s+N,t}} = 
\begin{cases}
\sum\limits_{l = 1}^{N} (G_{nls} + G_{nsl} ) V_{l+N,m}, & t = m,\ s =
1,\dots,N, \\ 
0, & t\neq m,%
\end{cases}%
\end{equation*}

The derivative of $\mathcal{R}(Q_h)$ can easily be calculated from \eqref{eq:reg}. 
\bigskip

\subsection{Algorithm}

The reconstruction of the unknown coefficient $c(x), \ x\in (0,b)$, is done
as follows.

\begin{itemize}
\item Step 1: Compute the Neumann data $g_1(k) = u_x(0,k)$ using \eqref{1.6}
then calculate $v_0(k) = \frac{g_1(k)}{k^2 g(k)}$.

\item Step 2: Compute $V^0$, $V^b$, and $\hat V_h$.

\item Step 3: Compute $Q_h$ by minimizing the cost function \eqref{eq:numqn}
and compute $V_h = Q_h + \hat V_h$.

\item Step 4: Compute $c(x), \ x\in (0,b),$ using \eqref{eq:c2} at $k = 
\underline{k}$.
\end{itemize}

Although the gradient projection algorithm is globally convergent, its
convergence is slow. Therefore, we use the Quasi-Newton method for
minimizing the objective functional.

\section{Numerical results}

\label{sec:result} In this section, we analyze the performance of the
proposed algorithm. For testing the algorithm against simulated data, we
solve the forward problem (\ref{1.1})--(\ref{1.2}) by converting it into the
1-D Lippmann-Schwinger equation 
\begin{equation*}
u(x,k)=u^{i}(x,k;x^{0})+k^{2}\int_{0}^{b}u^{i}(\xi ,k;x)[c(\xi )-1]u(\xi
,k)d\xi ,
\end{equation*}%
where $u^{i}(x,k;x^{0})=\frac{1}{2ik}e^{-ik|x-x^{0}|}$ is the incident wave
generated by the point source at $x=x^{0}$ in the homogeneous medium. This
integral equation is easily solved by approximating the integral by a
discrete sum.

\begin{figure}[tph]
\centering
\begin{tabular}{cc}
\includegraphics[width=0.5\textwidth]{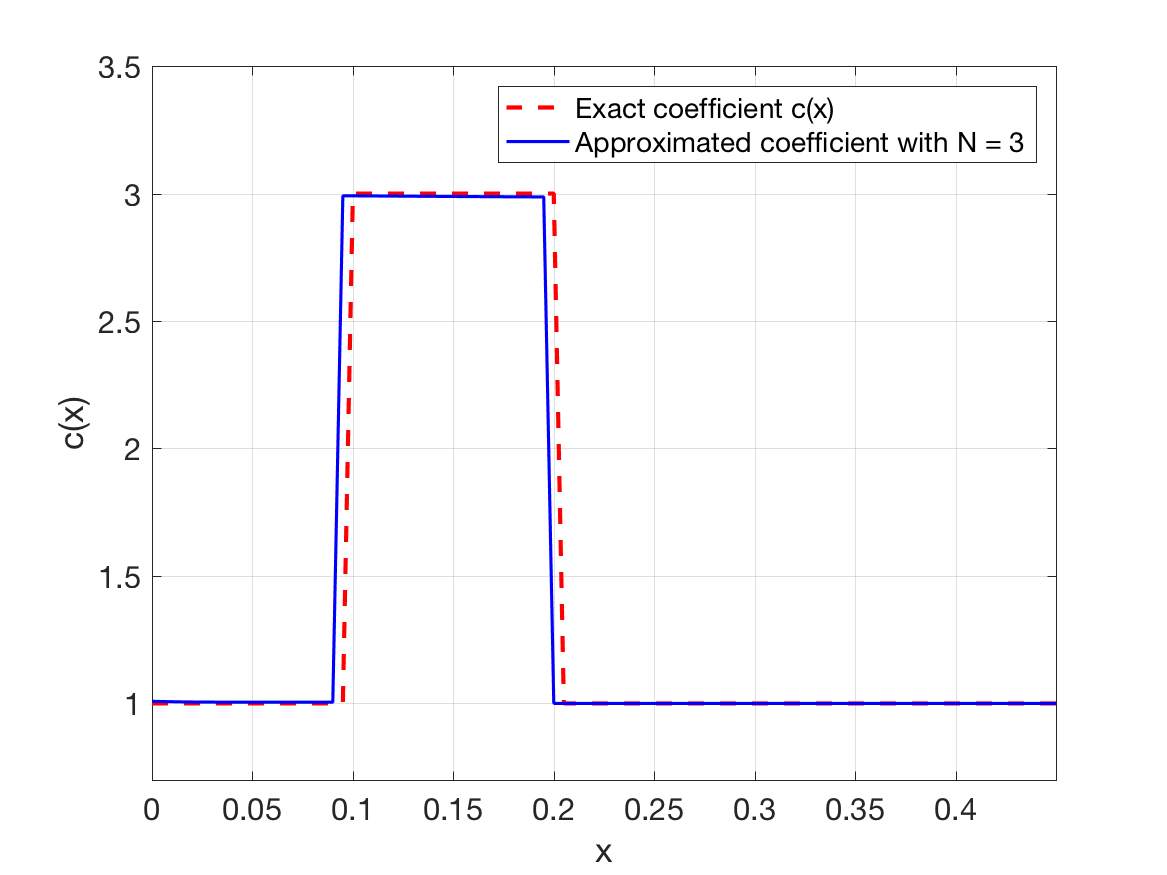}
& \includegraphics[width=0.5%
\textwidth]{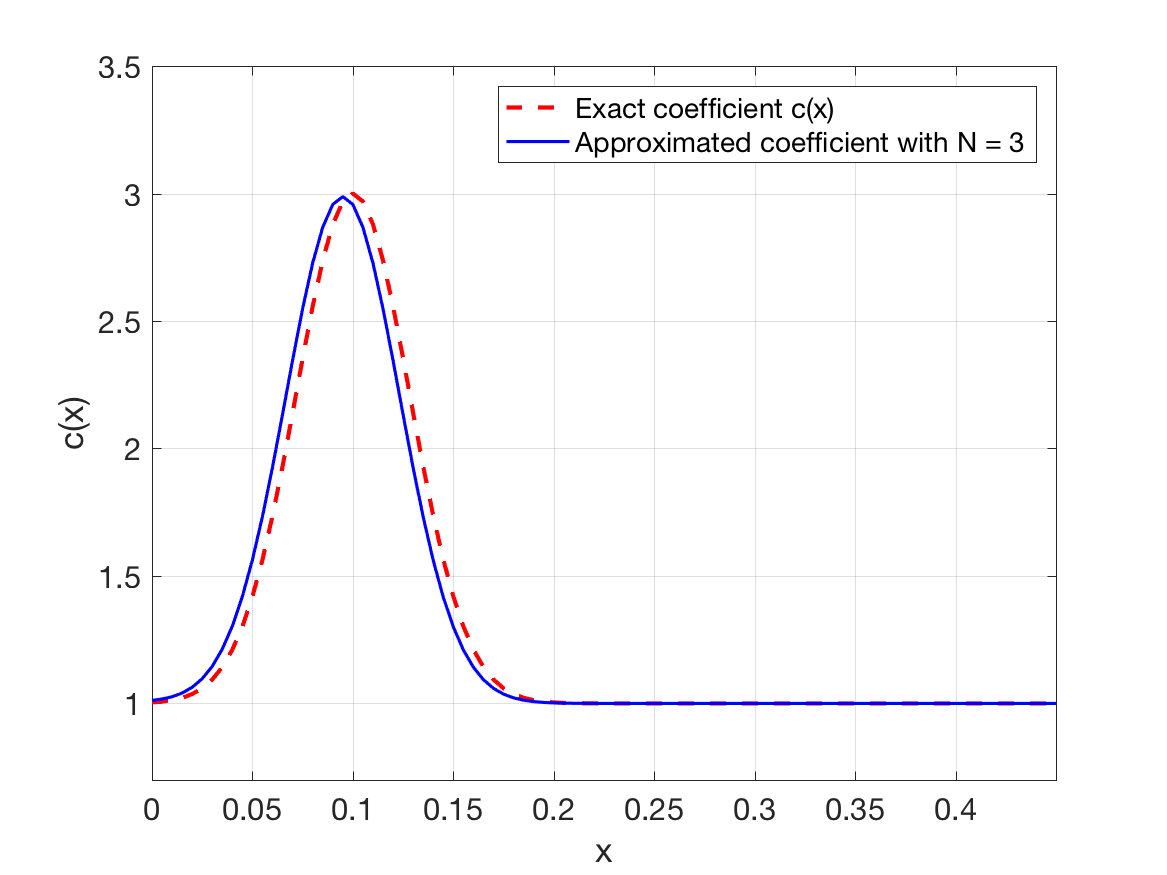} \\ 
a) $c(x)=1+3\chi[0.1,0.2],$ & b) $c(x)=1+2e^{-(x-0.1)^{2}/\left( 0.04\right)
^{2}}$%
\end{tabular}%
\caption{Comparison between the exact coefficient $c(x)$ and the approximate
coefficient computed by \eqref{eq:c2} with 3 basis functions and exact
functions $v_n$. The derivative $v^{\prime }_n(x)$ is approximated by a
finite difference quotient.}
\label{fig:1}
\end{figure}

In the following examples, to obtain the simulated data we solved the forward problem for $x\in [0,0.5]$ and added a $5\%$ of additive noise to the solution of the forward problem. That means, the noisy data $g_{noisy}(k)$ is calculated as $$g_{noise}(k)=g_{exact}(k) + 0.05 |g_{exact}(k) |\cdot \texttt{rand},$$
where $\texttt{rand}\in \left( -1,1\right) $  is a random
variable.

\begin{figure}[tph]
\centering
\includegraphics[width=0.5\textwidth]{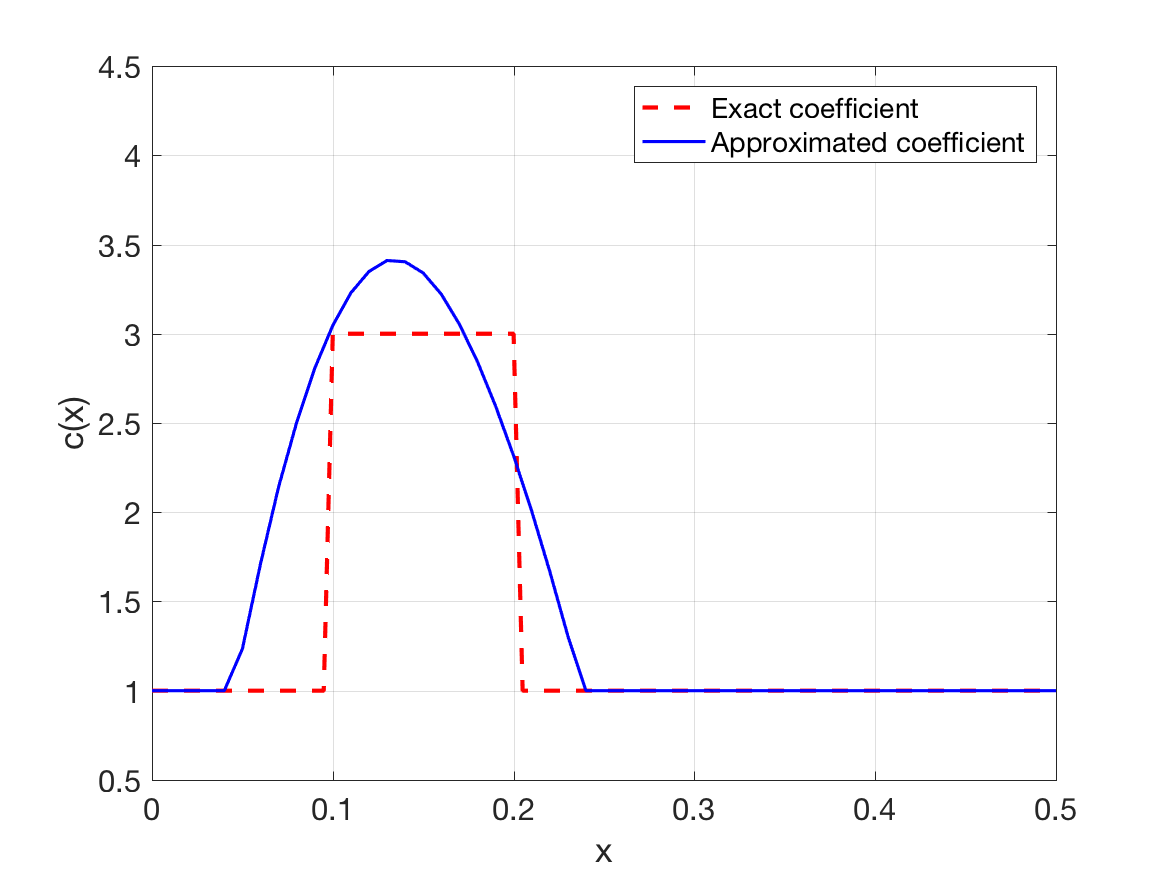}
\caption{Comparison of the exact and reconstructed coefficient $c(x)$ for
Example 1 for noiseless data and with the initial guess $Q_h^{(0)} = Q^*$.}
\label{fig:exa4-1}
\end{figure}

In solving the inverse problem, we chose the parameters as follows. Assume that $c(x)$ was unknown on $\left[ 0,b\right] =[0,0.3]$ only and $%
c(x)=1$ for $x\notin \lbrack 0,0.3]$, see (\ref{1.3}). The interval $[0,0.3]$ was divided into 31 subintervals of equal width $h = 0.01$. We
chose $11$ wavenumbers uniformly distributed between $\underline{k}=1$ and $%
\bar{k}=3$. Using numerical tests, we have observed that the coefficient $%
c(x)$ could be approximated quite accurately using only three terms in the
truncated Fourier series (\ref{2.10}) (see Figure \ref{fig:1}). Therefore,
the number of basis functions was chosen as $N=3$. The Carleman weight
coefficient was chosen as $\lambda =1$ and the regularization parameter was
chosen to be $\alpha =10^{-4}$. These parameters were chosen by
trial-and-error for a simulated data set. To analyze the reliability of the
algorithm, the same parameters were used for all other tests. 
Finally, the algorithm was started from the initial guess $Q^{(0)}_h\equiv 0$ in all of the following examples, except in Figure~\ref{fig:exa4-1} in which we show the effect of the truncation \eqref{2.10} on the reconstruction accuracy of $c(x)$.

Since we assume that $c(x) \ge 1$ for all $x$, we also replaced values of $%
c(x)$ which are less than 1 by 1. Note that this truncation was done as a
post-processing step after the objective functional $J_{h,\lambda,\alpha}$
was minimized. Therefore, it does not affect the inverse algorithm.

\bigskip

\begin{figure}[tph]
\centering
\includegraphics[width=0.5\textwidth]{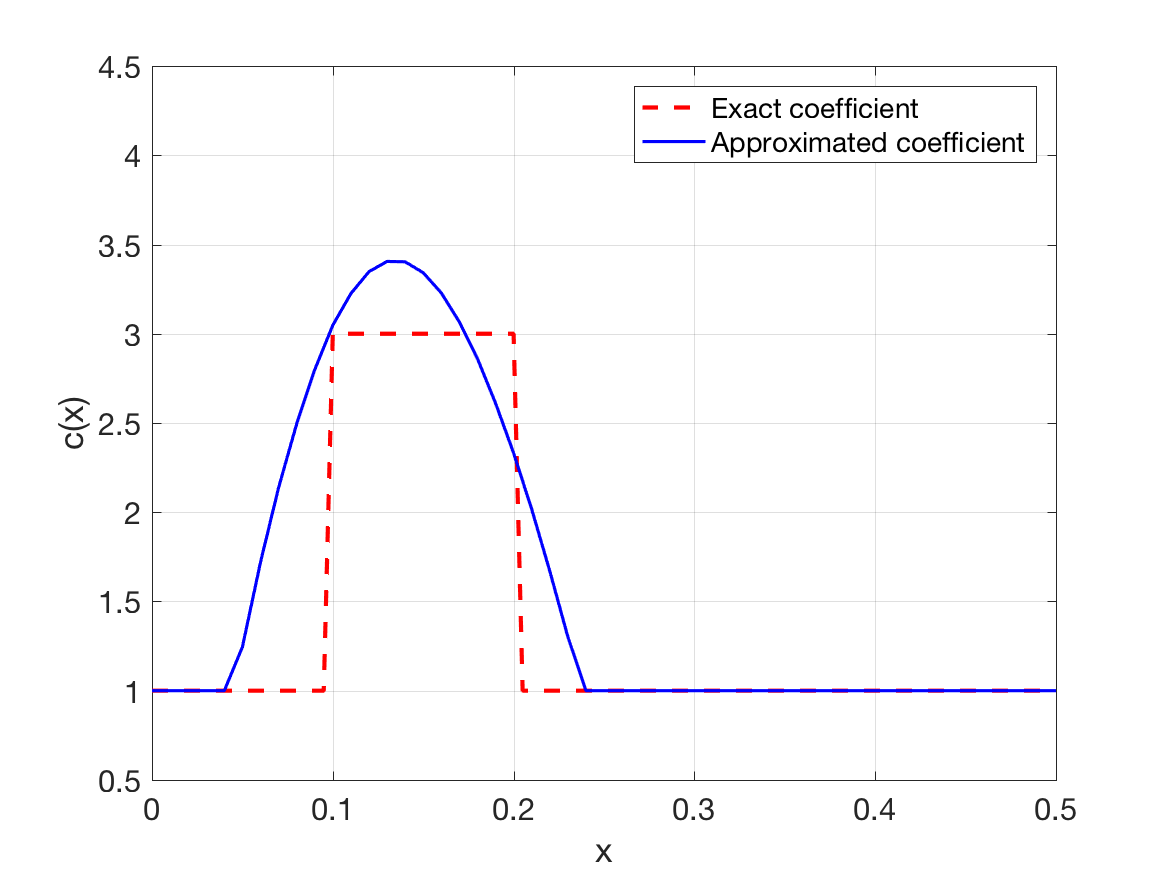}
\caption{Comparison of the exact and reconstructed coefficient $c(x)$ for
Example 1 for $5\%$ noisy data and the initial guess $Q_h^{(0)} \equiv 0$.}
\label{fig:exa1-1}
\end{figure}

\begin{figure}[tph]
\centering
\begin{tabular}{cc}
\includegraphics[width=0.5\textwidth]{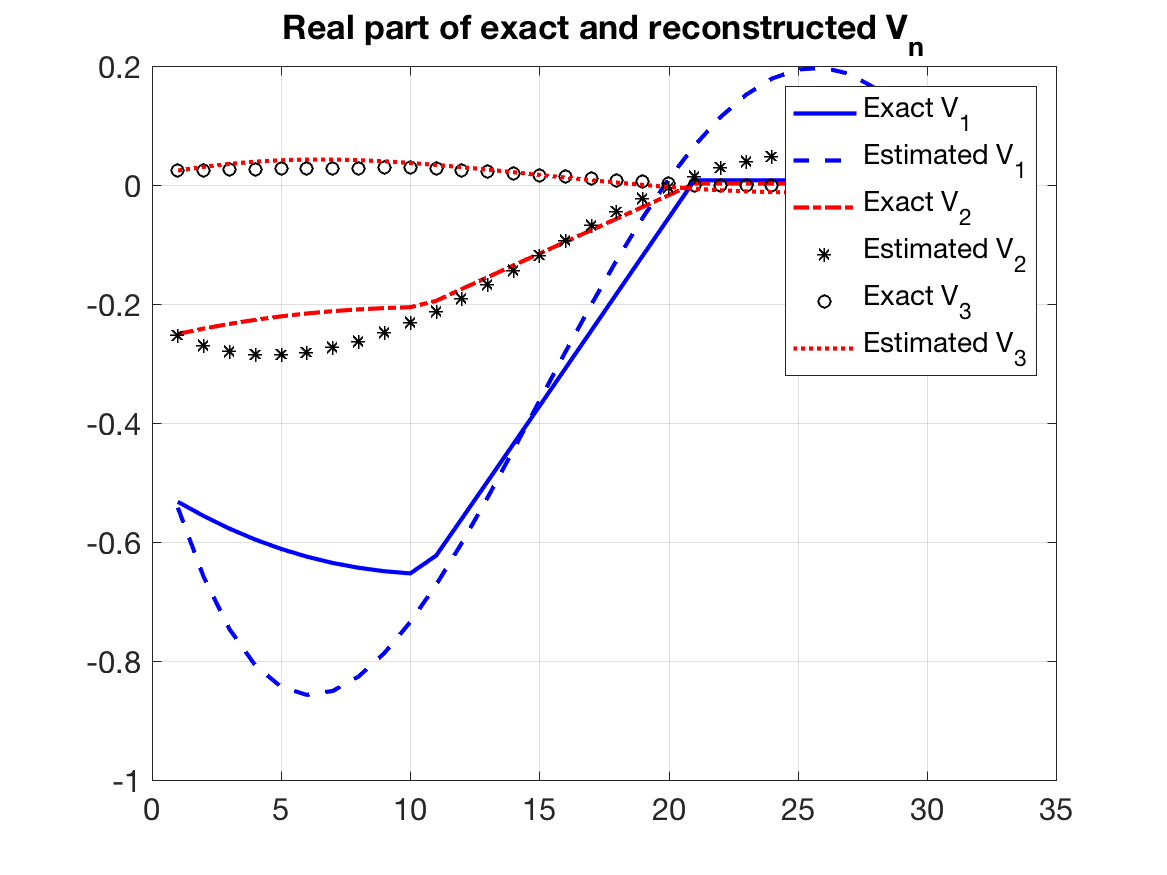} & %
\includegraphics[width=0.5\textwidth]{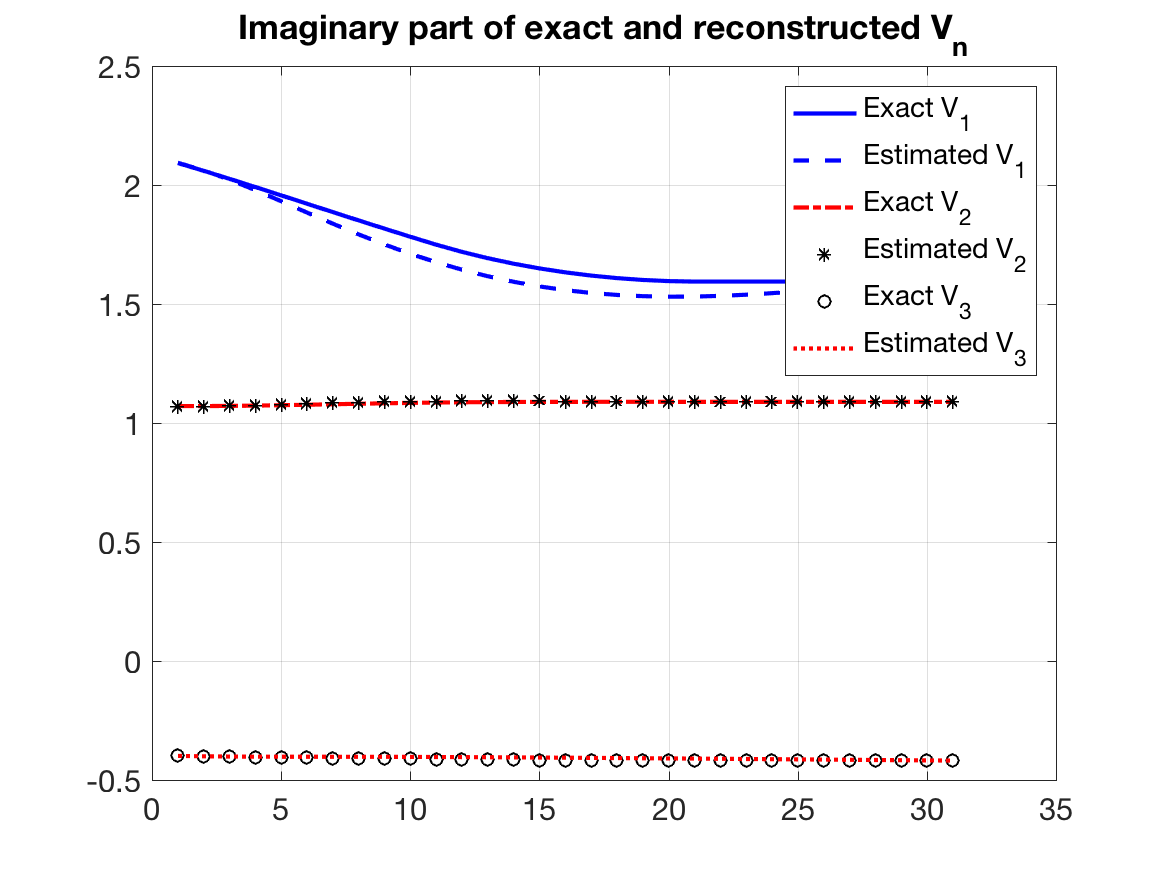} \\ 
a) Real part & b) Imaginary part%
\end{tabular}%
\caption{Comparison of the exact and reconstructed functions $V_n$ for
Example 1 for $5\%$ noisy data and the initial guess $Q_h^{(0)} \equiv 0$.}
\label{fig:exa1-2}
\end{figure}

\noindent \textbf{Example 1.} In the first example, we reconstruct a
piecewise constant coefficient given by $c(x)=1+3\chi \lbrack 0.1,0.2],$
where $\chi $ is the characteristic function. The reconstructed coefficient
is shown in Figure \ref{fig:exa1-1} together with the exact coefficient. In
Figure \ref{fig:exa1-2}, we show the reconstructed functions $V_{n}$, $%
n=1,2,3$, together with the \textquotedblleft exact\textquotedblright\ ones.
The exact functions $V_{n}$ are calculated from the solution of the forward
problem with the exact coefficient using \eqref{def:v}.

\begin{figure}[tph]
\centering
\begin{tabular}{cc}
\includegraphics[width=0.45\textwidth]{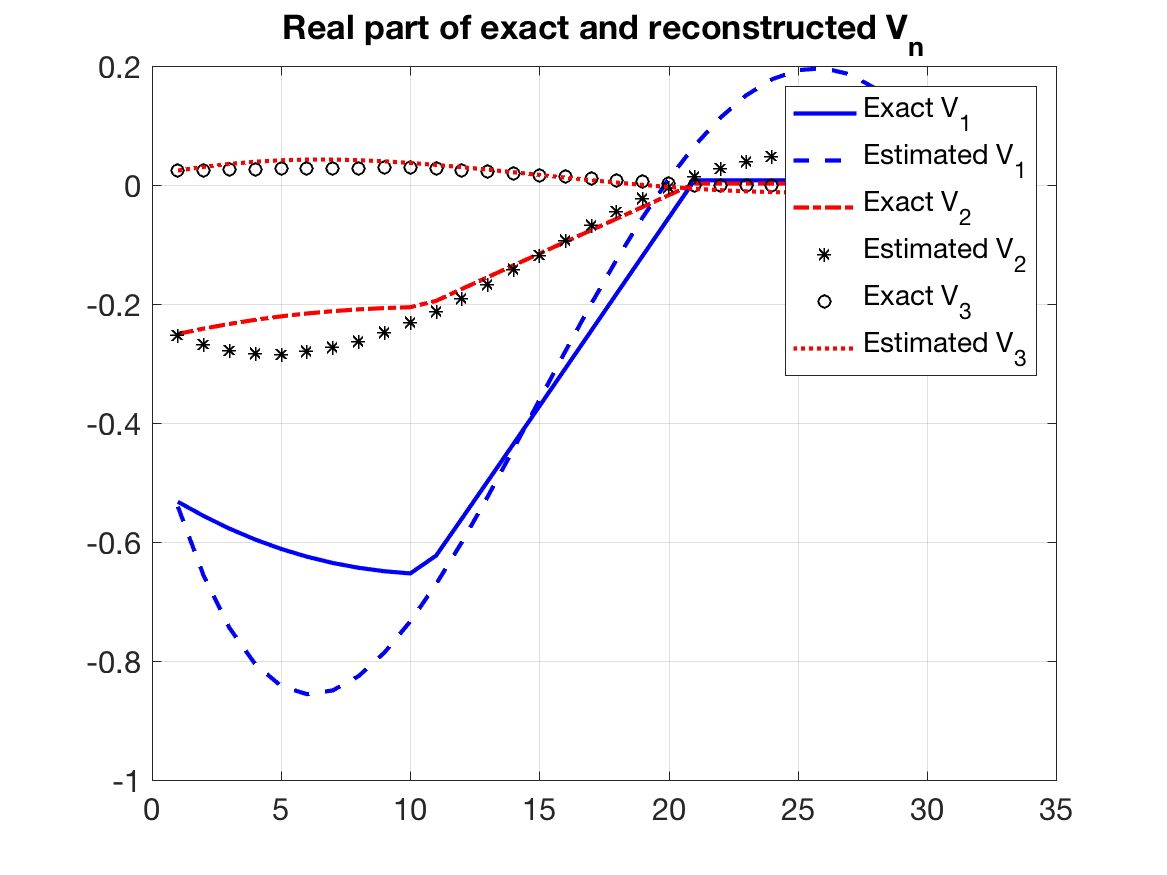} & %
\includegraphics[width=0.45\textwidth]{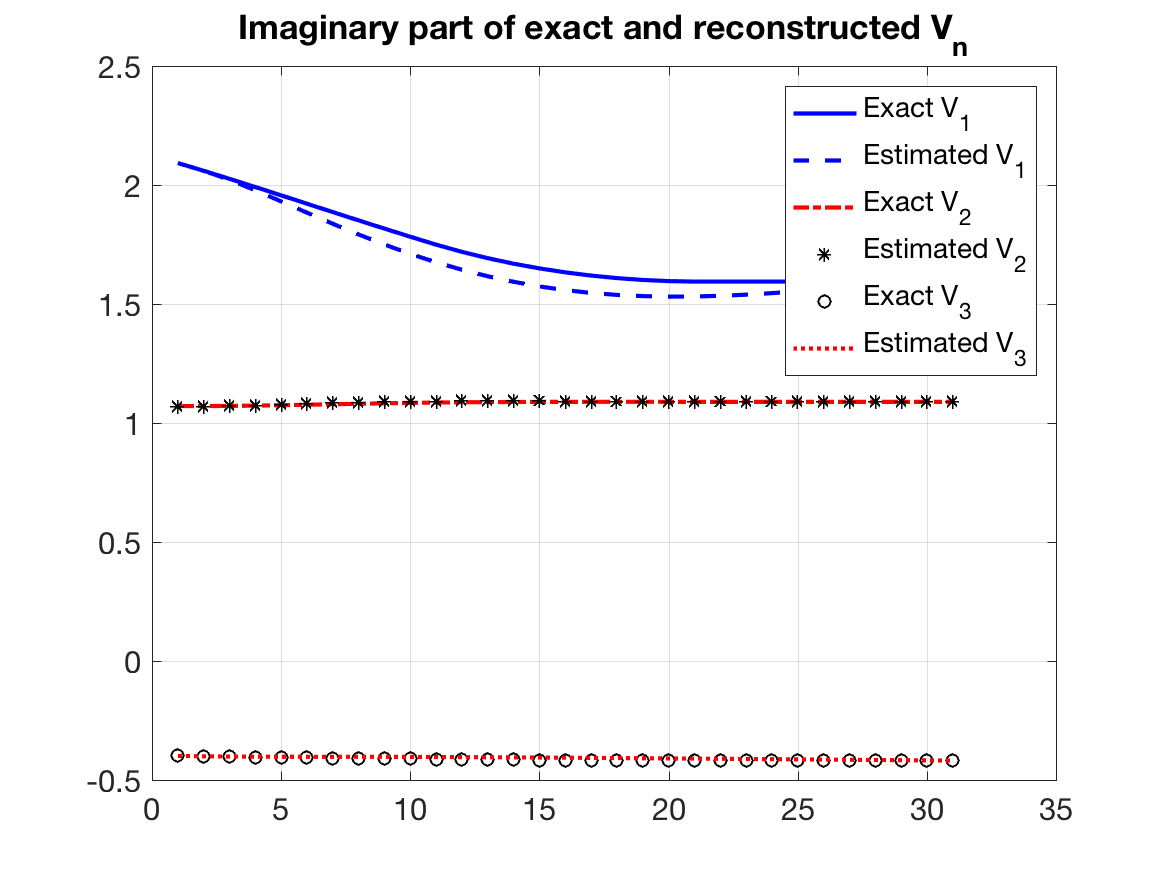} \\ 
a) Real part & b) Imaginary part%
\end{tabular}%
\caption{Comparison of the exact and reconstructed functions $V_n$ for
Example 1 for noiseless data and with the initial guess $Q_h^{(0)} = Q^*$.}
\label{fig:exa4-2}
\end{figure}

As can be seen from Figure \ref{fig:exa1-1}, the reconstructed coefficient
is a reasonable approximation of the true coefficient. One reason of the
difference between the exact and reconstructed coefficients that we have
observed through numerical analysis is due to the fact that \eqref{2.12} is
only an approximation. As a result, the \textquotedblleft
exact\textquotedblright\ functions $V_{n}$ are generally not the global
minimizer of $J_{h,\lambda ,\alpha }$. Therefore, when we minimize $%
J_{h,\lambda ,\alpha }$, we only obtain an approximation of these functions. To confirm this analysis, we show in Figures \ref{fig:exa4-1} and \ref{fig:exa4-2} the reconstructed coefficient $c(x)$ and the functions $V_n$ for noiseless data and with the initial guess $Q^{(0)}_h $ calculated from the exact function $V^*$. That means, $Q^{(0)}_h  = V^* - \hat V^*$. Figure \ref{fig:exa4-2} indicates that the exact $V_n$ are not the globally minimizer of $J_{h,\lambda ,\alpha }$ even with noiseless data. Comparison of Figures \ref{fig:exa1-1}, \ref{fig:exa1-2} with Figures \ref{fig:exa4-1}, \ref{fig:exa4-2} also  reveals another interesting observation that solutions resulted from the two initial guesses practically coincide.  This is exactly the thing which follows from our above theory.


One simple way to improve the
accuracy is to combine this globally convergent algorithm with a locally
convergent algorithm, such as the least-squares method. More precisely, we
can use the result of this algorithm as an initial guess for the
least-squares method. Since we want to focus on the performance of the
globally convergent algorithm, we do not discuss the least-squares method
here. We refer the reader to \cite{KT:SIAP2015} for this topic for a similar
problem in time domain.

%

\begin{figure}[tph]
\centering
\begin{tabular}{cc}
\includegraphics[width=0.5\textwidth]{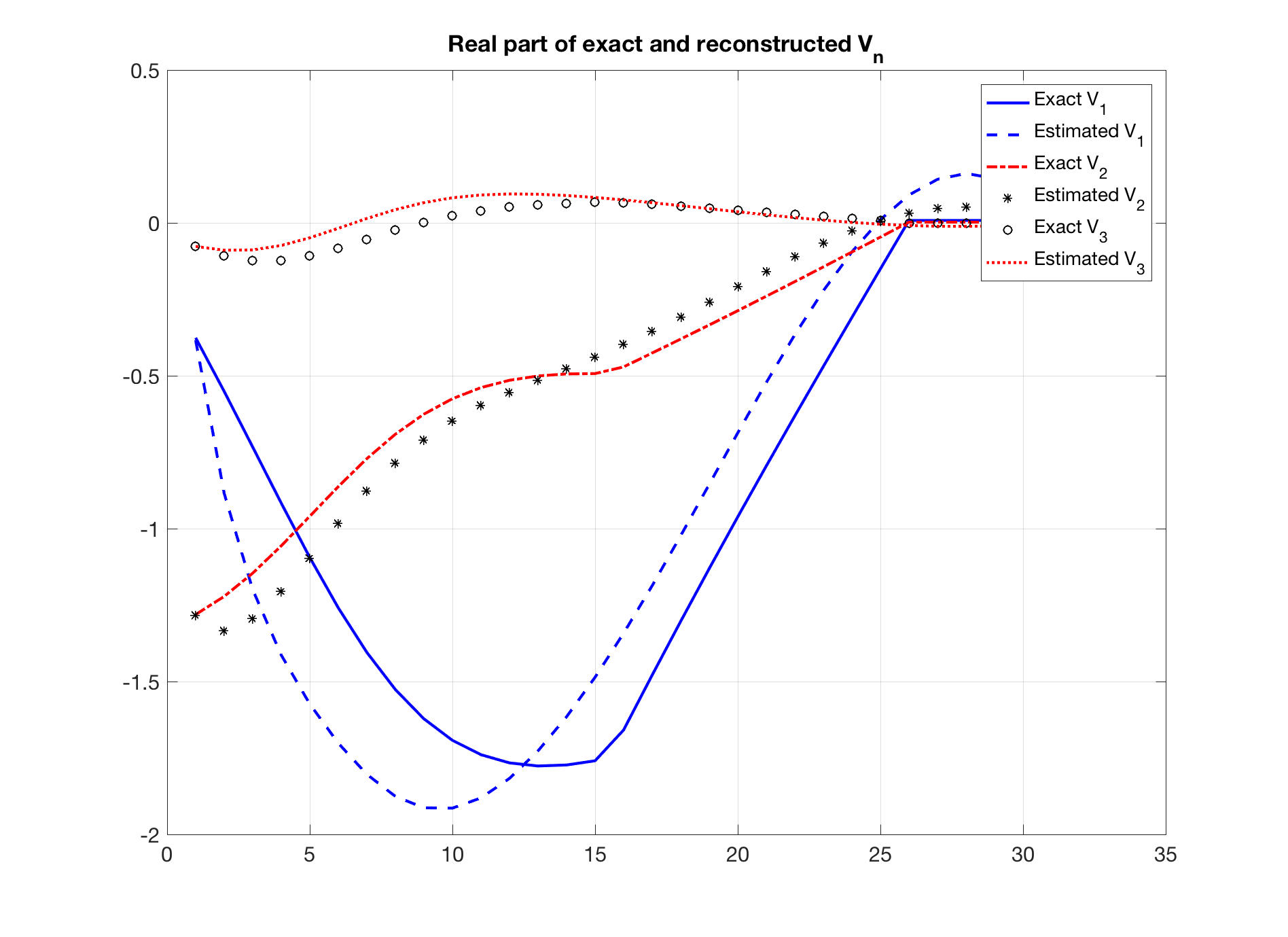} & %
\includegraphics[width=0.5\textwidth]{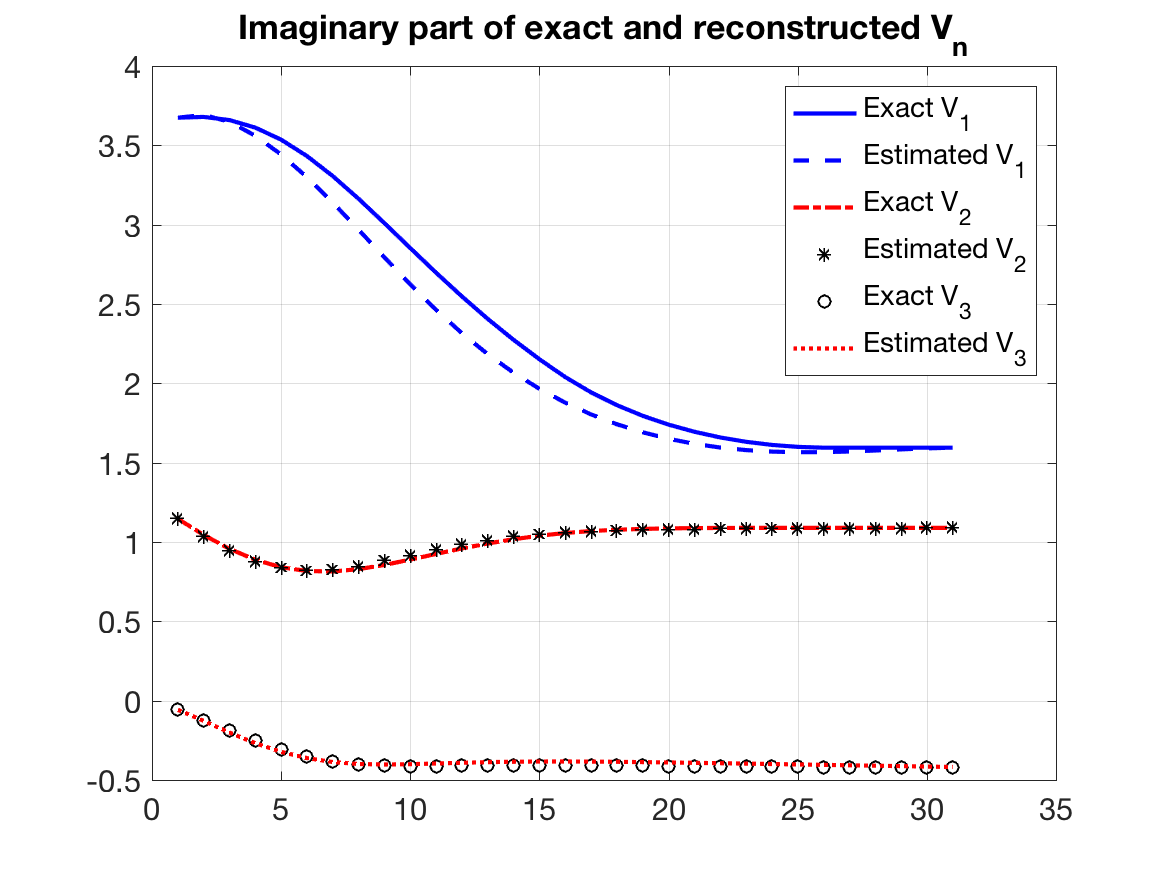} \\ 
a) Real part & b) Imaginary part%
\end{tabular}%
\caption{Comparison of the exact and reconstructed functions $V_n$ for
Example 2.}
\label{fig:exa2-2}
\end{figure}

\noindent \textbf{Example 2.} In this example, we consider another piecewise
constant coefficient with a larger inclusion/background contrast, $%
c(x)=1+6\chi \lbrack 0.15,0.25].$ The result is shown in Figures \ref%
{fig:exa2-1} and \ref{fig:exa2-2}. Even though the jump of the coefficient
is high in this case, we still can obtain the contrast quite well.

\begin{figure}[tph]
\centering
\includegraphics[width=0.5\textwidth]{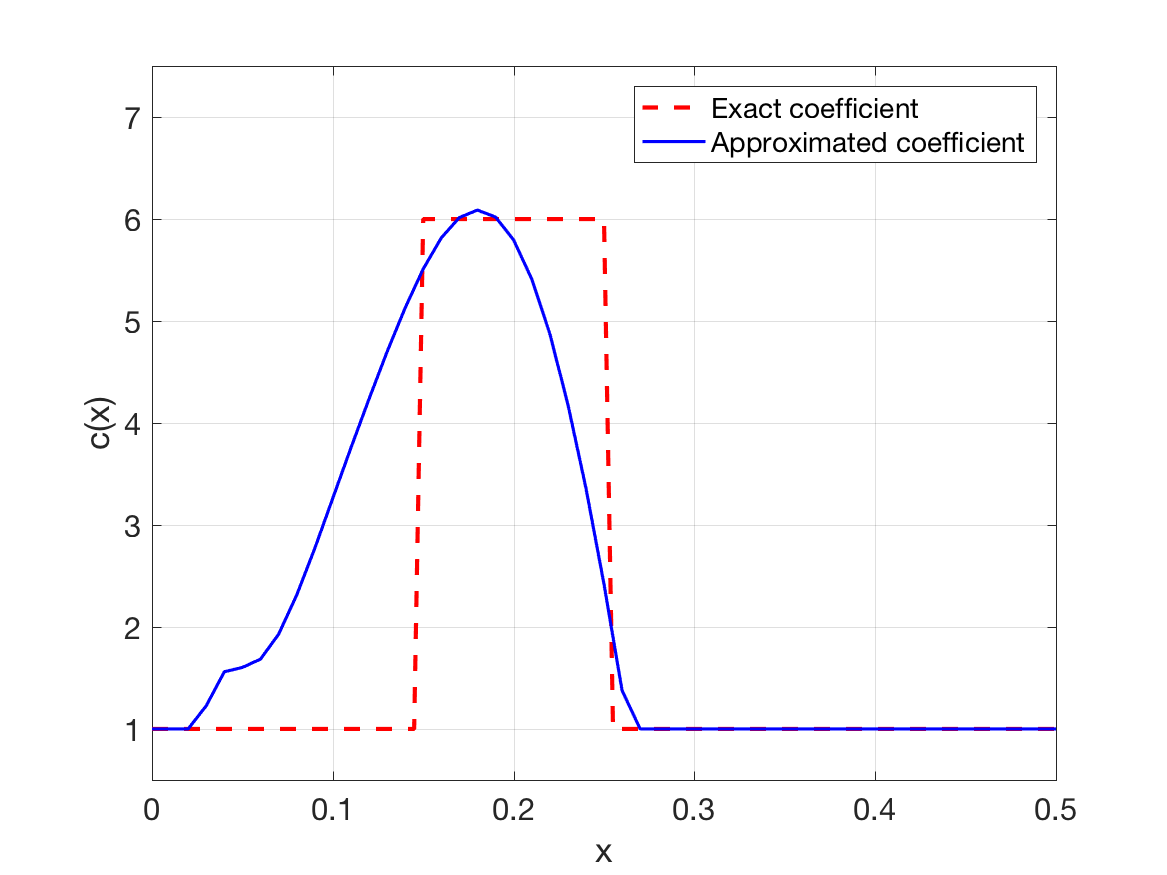}
\caption{Comparison of the exact and reconstructed coefficient $c(x)$ for
Example 2.}
\label{fig:exa2-1}
\end{figure}

\bigskip

\begin{figure}[tph]
\centering
\begin{tabular}{cc}
\includegraphics[width=0.5\textwidth]{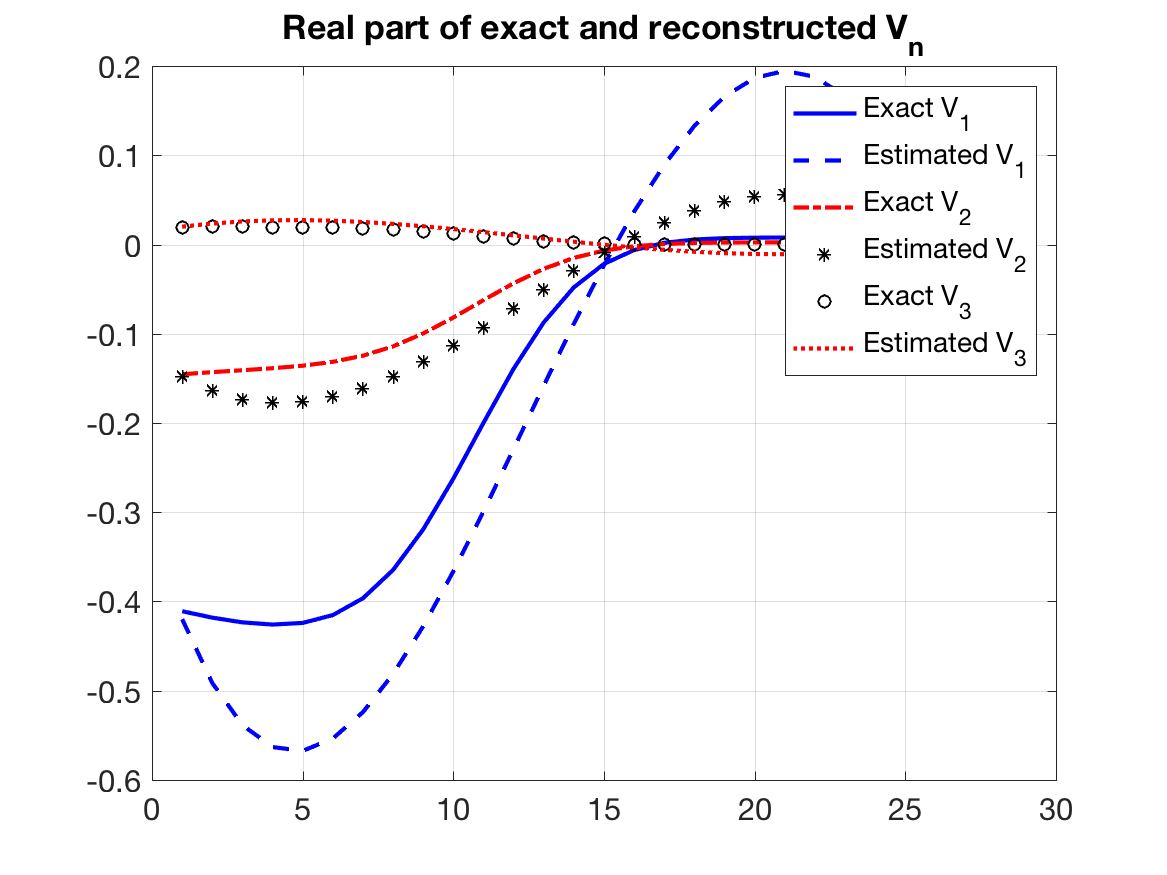} & %
\includegraphics[width=0.5\textwidth]{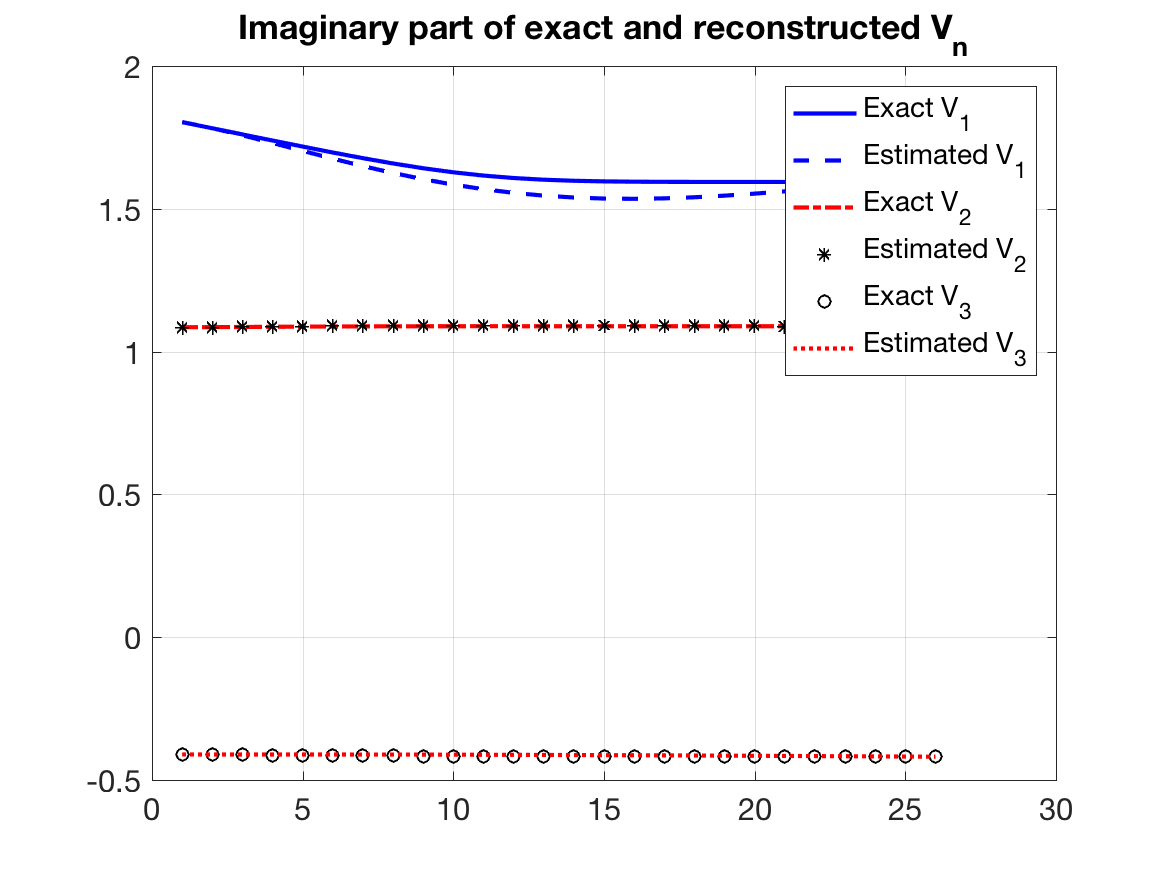} \\ 
a) Real part & b) Imaginary part%
\end{tabular}%
\caption{Comparison of the exact and reconstructed functions $V_n$ for
Example 3.}
\label{fig:exa3-2}
\end{figure}

\noindent \textbf{Example 3.} Finally, we consider a continuous coefficient
given by $c(x)=1+3e^{-(x-0.1)^{2}/\left( 0.04\right) ^{2}}.$ Figures \ref%
{fig:exa3-1} and \ref{fig:exa3-2} also show a reasonable reconstruction
result for both the coefficient $c(x)$ and the functions $V_n$.

\begin{figure}[tph]
\centering
\includegraphics[width=0.5\textwidth]{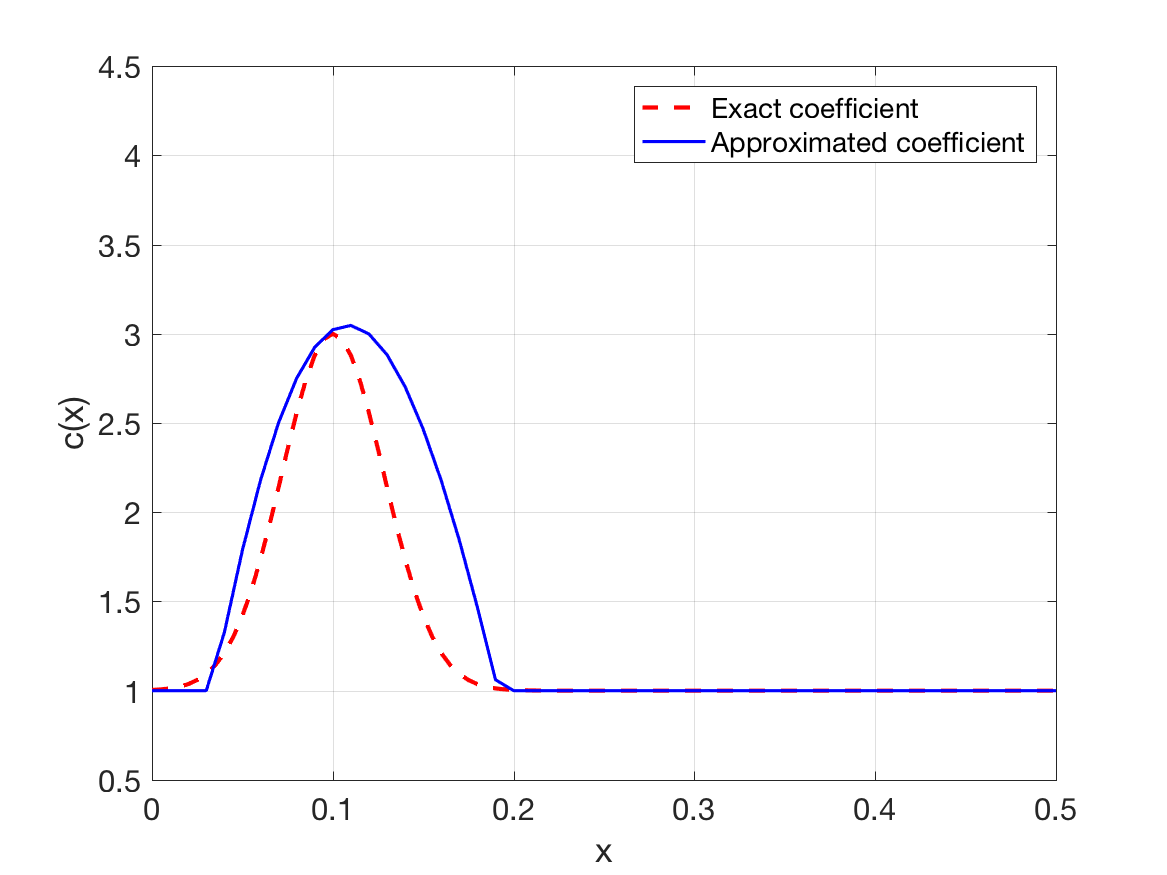}
\caption{Comparison of the exact and reconstructed coefficient $c(x)$ for
Example 3.}
\label{fig:exa3-1}
\end{figure}

\section{Concluding Remarks}

\label{sec:7} We have proposed a new globally convergent algorithm for the
multi-frequency inverse medium scattering problem. The main advantage of
this method is that we do not need a good first guess. The numerical
examples confirmed that the proposed method provides reasonable
reconstruction results. They also confirmed the global convergence of the proposed algorithm because the solutions from different initial guesses are practically the same. As a direct extension of this work, we are
considering the 2-d and 3-d problems and will report these cases in our
future work.

 \section*{Acknowledgment}
 The work of MVK was supported by US Army Research Laboratory and US Army Research Office grant W911NF-19-1-0044.

%

\end{document}